	\documentstyle{article}
	\newcommand{\al}{\alpha}
	\newcommand{\del}{\delta}
	
	\newcommand{\lam}{\lambda}

	\newcommand{\om}{\omega}
	\newcommand{\Del}{{\mathit \Delta}}
	\newcommand{\Gam}{{\mathit \Gamma}}
	
	\newcommand{\PHI}{{\mathit \Phi}}

	\font\smbbb=msbm5 
	\font\bbb=msbm7 
	\font\BBB=msbm10 
	\newcommand{\NN}{{\mbox{\BBB{N}}}}
	\newcommand{\ZZ}{{\mbox{\BBB{Z}}}}

	\newcommand{\re}{{\mbox{\bbb{R}}}}
	\newcommand{\RE}{{\mbox{\BBB{R}}}}
	\newcommand{\PP}{{\mbox{\BBB{P}}}}
	\newcommand{\smco}{{\mbox{\smbbb{C}}}}
	\newcommand{\co}{{\mbox{\bbb{C}}}}
	\newcommand{\CO}{{\mbox{\BBB{C}}}}
	\newcommand{\CP}{\PP}
	\newcommand{\CM}{\CO^{\times}}
	\newcommand{\cm}{{\co^{\times}}}
	\font\frak=eufm10 at 11 pt
	\font\sfrak=eufm9 
	\newcommand{\g}{\mbox{\frak{g}}}
	\newcommand{\gb}{\mbox{\frak{b}}}
	\newcommand{\gh}{\mbox{\frak{h}}}
	\newcommand{\gn}{\mbox{\frak{n}}}
	
	\newcommand{\gt}{\mbox{\frak{t}}}
	\newcommand{\gsl}{\mbox{\frak{sl}}}
	\newcommand{\sg}{\mbox{\sfrak{g}}}
	\newcommand{\sgb}{\mbox{\sfrak{b}}}
	\newcommand{\A}{{\cal A}}
	\newcommand{\CC}{{\cal C}}
	\newcommand{\FF}{{\cal F}}
	\newcommand{\CF}{\CC^*(\FF)}
	\newcommand{\HH}{{\cal H}}
	\newcommand{\LL}{{\cal L}}
	\newcommand{\OO}{{\cal O}}
	\newcommand{\UU}{{\cal U}}

	\newcounter{sect}
	\newcounter{subsect}
	\newcommand{\sect}[1]{\vspace{2ex}
		\addtocounter{sect}{1}\setcounter{subsect}{0}
		\begin{flushleft}
		{{\large\bf \arabic{sect}. {#1}}}
		\end{flushleft}
		\setcounter{thm}{0}
		\setcounter{equation}{0}
		\def\theequation{\arabic{sect}.\arabic{equation}}
		\def\thefigure{\arabic{sect}.\arabic{figure}}}
	\newcommand{\subsect}[1]{\vspace{0.25ex}\addtocounter{subsect}{1}
		\begin{flushleft}
		{{\bf \arabic{sect}.\arabic{subsect} {#1}}}
		\end{flushleft}}
	\newtheorem{thm}{Theorem}[sect]
	\newtheorem{prop}[thm]{Proposition}
	\newtheorem{lemma}[thm]{Lemma}
	\newtheorem{cor}[thm]{Corollary}
	\newtheorem{assump}[thm]{Assumption}
	\newtheorem{defn}[thm]{Definition}
	\newtheorem{ex}[thm]{Example}
	\newtheorem{rmk}[thm]{Remark}
	\newcommand{\proof}[1]{\noindent {\em Proof.}$\quad$ {#1} $\hfill\Box$}
	\newcommand{\be}{\begin{equation}}
	\newcommand{\ee}{\end{equation}}
	\newcommand{\bea}{\begin{eqnarray}}
	\newcommand{\eea}{\end{eqnarray}}
	\newcommand{\nno}{\nonumber \\}
	\newcommand{\sep}[1]{\!\!\!\! &{#1}& \!\!\!\! }
	\newcommand{\eq}{\sep{=}}
	\newcommand{\vc}{\sep{ }}

	\newcommand{\bra}{\langle}
	\newcommand{\ket}{\rangle}

	\newcommand{\medwedge}{\mbox{\fontsize{12pt}{0pt}\selectfont $\wedge$}}
	\newcommand{\inv}[1]{\frac{1}{#1}}
	\newcommand{\hf}{{\textstyle \inv{2}}}
	\newcommand{\e}[1]{e^{{#1}}}
        \newcommand{\ii}{\sqrt{-1}}
	\newcommand{\ch}{{\,\mathrm{char}\,}}
	
	\newcommand{\supp}{{\,\mathrm{supp}\,}}
	\newcommand{\spanc}{{\,\mathrm{span}\,}_\co}
	
	\newcommand{\cpct}{_{\mathrm{c}}}
	\newcommand{\set}[2]{\{{#1}\,|\,{#2}\}}
	\newcommand{\dr}{d}
	\newcommand{\zero}{{\{0\}}}
	
	\newcommand{\indlim}{{\displaystyle \lim_{\longrightarrow}}}
	\newcommand{\pdb}{\bar{\partial}}
	\newcommand{\ddz}{\frac{\partial}{\partial z}}
	\newcommand{\pq}{^{pq}}
	\newcommand{\pqr}{\pq_r}
	\newcommand{\Gamc}{\Gam\cpct}

	\newcommand{\ka}{K\"ahler }
	\newcommand{\bb}{Bia\l ynicki-Birula }
	\newcommand{\FFt}{\tilde{\FF}}
	\newcommand{\iii}{{i_0,\cdots,i_q}}
	\newcommand{\UUU}{U_{i_0}\cap\cdots\cap U_{i_q}}
	\newcommand{\two}[4]{\left\{    \begin{array}{ll}
					{#1}, & {\mbox{if }} {#2}, \\
					{#3}, & {\mbox{if }} {#4}
					\end{array}     \right.}
	\newcommand{\three}[3]{\left\{	\begin{array}{l}
					{#1}	\\	{#2}	\\	{#3}
					\end{array}	\right.}

\begin{document}
$\!\,{}$

        \vspace{-5ex}

        \begin{flushright}
{\tt math.AG/9806118} (June, 1998)\\
Adelaide IGA preprint 1998-07
        \end{flushright}

\vspace{1.25ex}

	\begin{center}
{\Large\bf On the Instanton Complex of Holomorphic Morse Theory}\\

	\vspace{4ex}
	{\large\rm Siye Wu}\footnote{E-mail address:
	{\tt swu@maths.adelaide.edu.au}}

	{\em Department of Pure Mathematics, University of Adelaide, 
	Adelaide, SA 5005, Australia}

	\end{center}

	\vspace{3ex}

        \begin{quote}
{\small {\bf Abstract.}
Consider a holomorphic torus action on vector bundles over a complex manifold
which lifts to a holomorphic vector bundle.
When the connected components of the fixed-point set are partially ordered,
we construct, using sheaf-theoretical techniques, two spectral sequences
that converges to the twisted Dolbeault cohomology groups and those with
compact support, respectively.
These spectral sequences are the holomorphic counterparts of the instanton
complex in standard Morse theory.
The results proved imply holomorphic Morse inequalities and fixed-point
formulas on a possibly non-compact manifold.
Finally, a number of examples and applications are given.}

        \end{quote}

        \vspace{3ex}

\sect{Introduction}

Given a Morse function on a compact real manifold, the Morse inequalities
bound the Betti numbers in terms of the information of critical points.
However, the former can not be determined by the Morse inequalities alone
unless the Morse function is perfect.
If the Morse function satisfies the transversality condition [\ref{Sm}],
then there is a finite dimensional complex, called the Thom-Smale-Witten
complex or the instanton complex [\ref{W1}], which computes 
the cohomology groups of the manifold. 
(See [\ref{B2}] for a historical review.)
The instanton complex consists of vector spaces spanned by the critical points
of the Morse function (when they are isolated), graded by their Morse indices.
The coboundary operators come from counting (with orientation) the number of 
gradient paths between critical points whose Morse indices differ by one.
The latter is related to the instanton tunneling effect in supersymmetric
quantum mechanics [\ref{W1}].

Consider a complex manifold with a holomorphic group action and a holomorphic
vector bundle over the manifold on which the group action lifts 
holomorphically.
We want to determine the Dolbeault cohomology groups (twisted by the vector 
bundle) as representations of the group.
When the manifold is compact, the fixed-point formula of Atiyah and Bott 
[\ref{AB}] (for isolated fixed points) and of Atiyah and Singer [\ref{AS}]
computes the alternating sum of the characters on the cohomology groups.
For holomorphic Morse theory, this (equivariant) index theorem is the
counterpart of the Hopf (or Lefschetz) formula.
When the manifold is compact and K\"ahler and the group is the circle group,
Morse-type inequalities were obtained by Witten [\ref{W2}] using a holomorphic
version of supersymmetric quantum mechanics.
These (equivariant) holomorphic Morse inequalities put constraints on 
the sizes of Dolbeault cohomology groups but do not completely determine them.
In [\ref{MW}], a heat kernel proof was given under the additional assumption
that the fixed points are isolated.
In [\ref{Wu}], these inequalities were generalized to cases with torus and
non-Abelian group actions.
Furthermore, it was shown that the \ka assumption was necessary for holomorphic
Morse inequalities [\ref{Wu}], although not so for the fixed-point theorem.
In [\ref{WZ}], these inequalities were proved analytically for compact \ka 
manifolds with possibly non-isolated fixed points.

In this paper, we construct the holomorphic counterpart of the instanton 
complex which computes the Dolbeault cohomology groups using 
the combinatorial data of the group action.
At the same time, we investigate more closely the condition on the 
complex manifold for establishing a holomorphic Morse theory.
Holomorphic Morse theory differs from ordinary Morse theory in a number of 
ways.
If the circle group acts on a compact \ka manifold in a Hamiltonian fashion,
the moment map is a perfect Morse function whose critical points have even
Morse indices only, which can not differ by one.
Furthermore, Smale's transversality condition fails in general and the
gradient paths are never isolated because of the circular symmetry.
Consequently, the techniques for holomorphic Morse theory will be
quite different from those for ordinary Morse theory.

We start with a complex manifold with a holomorphic action of a (complex)
torus.
The action of a non-compact $1$-parameter subgroup is analogous to the
gradient flow of a Morse function.
The group action is meromorphic if, roughly speaking, all such orbits start
from and end at some points in the manifold, which must be fixed points of
the torus.
If so, then there is a relation on the connected components of the 
fixed-point set given by the direction of the flows.
The central result of this paper is that if this relation is a partial
ordering, then there are two (equivariant) spectral sequences converging
(equivariantly) to the twisted Dolbeault cohomology groups and those with
compact support, respectively.
These spectral sequences will be constructed using sheaf-theoretic techniques
from a filtration of the complex manifold determined by the group action.
The spectral sequences, with the natural coboundary maps, are the couterparts
in holomorphic Morse theory of the instanton complex in ordinary Morse theory.
The information of the $E_1$-terms already implies the holomorphic Morse 
inequalities.
But unlike ordinary Morse theory, the spectral sequences do not always
degenerate at $E_2$.
When the manifold is compact and K\"ahler, the partial order condition 
is automatically satisfied.
Thus the results of [\ref{W2}, \ref{MW}, \ref{Wu}, \ref{WZ}] are recovered.

The rest of the paper is organized as follows.
In section~2, we establish some facts about meromorphic torus actions
on a compact or a suitably non-compact complex manifold.
In section~3, we construct two spectral sequences converging to 
Dolbeault cohomology groups and those with compact support, respectively,
under the partial order condition. 
In particular, we obtain holomorphic Morse inequalities and fixed-point
formulas for a possibly non-compact manifold.
We also study the condition under which the spectral sequences degenerate
to cochain complexes.
In section~4, we consider various examples and applications.
We first present a spectral sequence calculation using the language
of \v Cech cohomology.
The application to flag manifolds yields a geometric realization
of the Bernstein-Gelfand-Gelfand resolution and its generalizations.
We also study the Dolbeault cohomologies and geometric quantization on
non-compact manifolds.

Throughout this paper, $\NN$, $\RE$, $\RE^\pm$, $\CO$ and $\CO^\times$
denote the sets of non-negative integers, real numbers, positive (negative)
real numbers, complex numbers and non-zero complex numbers, respectively.

\sect{Holomorphic torus actions}

We first recall from [\ref{Wu}, \ref{WZ}] some notations 
of holomorphic torus actions without making the compact or \ka assumption.

Let $T$ be a complex torus with Lie algebra $\gt$.
Let $T_\re$ be the (real) maximal compact torus subgroup of $T$
and $\gt_\re=\ii\,{\rm Lie}(T_\re)$.
Let $\ell$ be the integral lattice in $\gt_\re$, and $\ell^*\subset\gt^*_\re$,
the dual lattice.
If $T=\CM$, the multiplicative group of non-zero complex numbers,
then $T_\re=S^1$, $t_\re=\RE$, and $\ell=\ZZ$.
In general, for any $v\in\ell-\zero$, there is an embedding $j_v\colon\CM\to T$
whose image $\CM_v$ is a $\CM$-subgroup of $T$.

The ring of formal characters of $T$ is
$\ZZ[\ell^*]=\set{q=\sum_{\xi\in\ell^*}q_\xi\e{\xi}}{q_\xi\in\ZZ}$.
The {\em support} of $q\in\ZZ[\ell^*]$ is
$\supp q=\set{\xi\in\ell^*}{q_\xi\ne0}$.
We say that $q\ge0$ if $q_\xi\ge0$ for all $\xi\in\ell^*$.
Consider a representation $R$ of $T$.
If every weight $\xi\in\ell^*$ of $R$ has a finite multiplicity $r_\xi$,
then the character $\ch R=\sum_{\xi\in\ell^*}r_\xi\e{\xi}\in\ZZ[\ell^*]$
is well-defined.
Let $\supp R=\supp\ch R$.
As in [\ref{MW}, \ref{Wu}], we write
	\be
\frac{\e{\eta}}{1-\e{\xi}}\stackrel{{\rm def.}}{=}
\sum_{k=0}^\infty\e{\eta+k\xi},\quad\quad\xi,\eta\in\ell^*.
	\ee
We emphasize here that the left-hand side is a notation for the formal series
in $\ZZ[\ell^*]$ on the right-hand side.
More generally, if $R$ is a finite dimensional representation of $T$,
we can write
	\be
\ch\inv{\det(1-R)}=\sum_{k=0}^\infty\ch S^k(R)=\ch S(R).
	\ee

Let $X$ be a complex manifold of dimension $n$.
Suppose $T$ acts holomorphically and effectively on $X$.
The fixed-point set $X^T$ of $T$ in $X$, if non-empty,
is a complex submanifold of $X$.
Let $F$ be the set of connected components of $X^T$.
Then $X^T=\bigcup_{\al\in F}X^T_\al$, where $X^T_\al$ is the component
labeled by $\al\in F$.
Let $n_\al=\dim_\co X^T_\al$.
Let $N_\al\to X^T_\al$ be the (holomorphic) normal bundle of $X^T_\al$ in $X$.
$T$ acts on $N_\al$ preserving the base $X^T_\al$ pointwise.
The weights of the isotropy representation on the normal fiber
remain constant within any connected component.
Let $\lam_{\al,k}\in\ell^*-\zero\subset\gt_\re^*$ ($1\le k\le n-n_\al$)
be the isotropy weights on $N_\al$.
The hyperplanes $(\lam_{\al,k})^\perp\subset\gt_\re$ cut $\gt_\re$ into open
polyhedral cones called {\em action chambers} [\ref{PW}].
Choose an action chamber $C$.
Let $\lam^C_{\al,k}=\pm\lam_{\al,k}$,
with the sign chosen so that $\lam^C_{\al,k}\in C^*$.
(Here $C^*$ is the dual cone in $\gt_\re^*$ defined by
$C^*=\set{\xi\in\gt_\re^*}{\bra\xi,C\ket>0}$.)
We define $\nu^C_\al$ as the number of weights $\lam_{\al,k}\in C^*$.
Let $N^C_\al$ be the direct sum of the sub-bundles corresponding to
the weights $\lam_{\al,k}\in C^*$.
Then $N_\al=N^C_\al\oplus N^{-C}_\al$.
$\nu^C_\al$ is the rank of the holomorphic vector bundle $N^C_\al$;
that of $N^{-C}_\al$ is $\nu^{-C}_\al=n-n_\al-\nu^C_\al$,
which is called the {\em polarizing index} of $X^T_\al$ with respect to $C$.

In subsection~\arabic{sect}.1, we will consider holomorphic torus actions
on compact manifolds;
a non-compact setting will be studied in subsection~\arabic{sect}.2.

\subsect{Meromorphic torus actions on compact manifolds}

Throughout this subsection, $X$ is a compact complex manifold with 
a holomorphic action of the torus $T$.
Then $F$ is a finite set and each component $X^T_\al$ ($\al\in F$) is compact.

\begin{defn}
{\em A holomorphic $T$-action on $X$ is {\em meromorphic} if for any $x\in X$
and any $v\in\ell-\zero$, the limit $\pi^v(x)=\lim_{u\to0}j_v(u)x$ exists.}
\end{defn}

If $T=\CM$, the action is meromorphic if and only if for any $x\in X$,
the limits $\pi^+(x)=\lim_{u\to0}ux$ $\pi^-(x)=\lim_{u\to\infty}$ exist.
In this case, the holomorphic map $\CM\times X\to X$ can be extended to 
a meromorphic map $\CP^1\times X\to X$.

\begin{prop}\label{LIMIT}
If the $T$-action on $X$ is meromorphic, then\\
1. for any $v\in\ell-\bigcup_{\al\in F,1\le k\le n-n_\al}(\lam_{\al,k})^\perp$,
the fixed-point set of $\CM_v$ coincides with $X^T$;\\
2. for any $x\in X$ and action chamber $C$, the limit $\pi^v(x)$ for 
$v\in\ell\cap C$ depends only on $C$ and not on the choice of $v$.
\end{prop}

\proof{1. Let $X'$ be a connected component of the fixed-point set of $\CM_v$.
Then $X'\cap X^T$ is a closed subset of $X'$.
For any $x\in X'\cap X^T$, let $X^T_\al$ be the component of $X^T$ 
that contains $x$.
Since the $T$-action is effective and $\lam_{\al,k}(v)\ne0$ for any 
$0\le k\le n-n_\al$, we have $\dim X^T_\al=\dim X'$.
Therefore $X'\cap X^T$ is also an open subset of $X'$.
Finally, choose $v_1,\dots,v_{r-1}\in\ell$ ($r=\dim_\co T$) such that
$\{v,v_1,\dots,v_{r-1}\}$ is a basis of $\gt_\re$.
Pick any $x'\in X'$.
Since the $T$-action is meromorphic, the iterated limit
$x=\pi^{v_1}\pi^{v_2}\cdots\pi^{v_{r-1}}(x')$ exists.
It is clear that $x\in X'\cap X^T$.
So $X'\cap X^T\ne\emptyset$.
Consequently, $X'\cap X^T=X'=X^T_\al$.\\
2. From part~1, we have $y=\pi^v(x)\in X^T_\al$ for some $\al\in F$.
By [\ref{CS1}, Proposition~I], there is a $T_\re$-invariant neighborhood
$W_y$ of $y$ in $N_\al$ and a $T$-equivariant holomorphic embedding
$\psi_y\colon W_y\to X$.
Let $X^v_y=(\pi^v)^{-1}(X^T_\al)$.
Then from the linear $T$-action on $N_\al$, we get 
$X^v_y\cap\psi_y(W_y)=\psi_y(N^C_\al\cap W_y)$.
Hence $X^v_y=T\,\psi_y(N^C_\al\cap W_y)$;
this depends only on $C$ and not on the choice of $v$.}

\vspace{1ex}
We denote $\pi^v(x)$ by $\pi^C(x)$ when $v\in\ell\cap C$.

\begin{rmk}{\em
1. If $X$ is a compact \ka manifold and $X^T\ne\emptyset$,
then the $T_\re$-action is Hamiltonian [\ref{F}].
Let $\mu\colon X\to\gt^*_\re$ be a moment map.
For $v\in\ell-\zero$, the $1$-parameter subgroup $\set{j_v(\e{t})}{t\in\RE}$
generates the gradient flows of $\bra\mu,v\ket$, along which its value
strictly decreases.
Therefore the limit $\pi^v(x)$ for any $x\in X$ exists and the $T$-action
is meromorphic.\\
2. A holomorphic action on $X$ may not be meromorphic even if $X$
is compact and K\"ahler.
For example, let $\ZZ$ act on $\CO-\zero$ by $k\colon z\mapsto 2^kz$
($k\in\ZZ$, $z\in\CO-\zero$) and let $X=(\CO-\zero)/\ZZ$ be the quotient.
Then the standard multiplication of $\CM$ on $\CO-\zero$ induces a holomorphic
action on $X$ which has no fixed points and hence is not meromorphic.}
\end{rmk}

In order to capture the topology of $X$ by the fixed-point information,
it is necessary to assume that the $T$-action is meromorphic.
If so, then $X$ has a cell decomposition according to the connected components
of $X^T$ that $\pi^C$ maps to.

\begin{defn}
{\em Suppose the $T$-action on $X$ is meromorphic.
Set $X^C_\al=(\pi^C)^{-1}(X^T_\al)$.
Then
	\be\label{bbc}
X=\bigcup_{\al\in F}X^C_\al
	\ee
is called the {\em \bb decomposition} with respect to $C$.}
\end{defn}

Consider the case $T=\CM$.
If the $\CM$-action is meromorphic, set $X^\pm_\al=(\pi^\pm)^{-1}(X^T_\al)$.
The decompositions $X=\bigcup_{\al\in F}X^\pm_\al$ are called the
{\em plus} ({\em minus}) decompositions, respectively.
For example, $X=\CP^1$ is the union of $\CO$ (with the standard
multiplication of $\CM$) and $\{\infty\}$. 
The fixed-point set is $X^T=\{0,\infty\}$.
We have $X^+_0=X-\{\infty\}$, $X^+_\infty=\{\infty\}$
and $X^-_0=\zero$, $X^-_\infty=X-\zero$.
Both decompositions $X=X^\pm_0\cup X^\pm_\infty$ consist of a $0$-cell
and a $2$-cell.

The cells $X^C_\al$ are $T$-invariant.
If the transversality condition is satisfied, then the decomposition 
(\ref{bbc}) is a stratification [\ref{BB2}, Theorem~5].
In general, this is not true even when $X$ is K\"ahler.
An example is the Hirzebruch surface (the blow-up of $\CP^2$ at one point)
[\ref{BB2}, Example 1].

\begin{defn}
{\em For $\al,\beta\in F$, we write $\al\to\beta$ if there is $x\in X$
such that $\pi^C(x)\in X^T_\al$ and $\pi^{-C}(x)\in X^T_\beta$.
We write $\al\prec\beta$ if either $\al=\beta$ or there is a {\em chain} 
from $\al$ to $\beta$,
i.e., a finite sequence $\al_0=\al, a_1,\dots,\al_{r-1},\al_r=\beta$
in $F$ such that $\al_{i-1}\to\al_i$ for all $1\le i\le r$ ($r>0$).
Such a chain is called a {\em quasicycle} of length $r$ if $\al=\beta$.}
\end{defn}

Obviously, the relation $\prec$ on $F$ depends on the choice of $C$.

\begin{rmk}{\em
Results on meromorphic $\CM$-actions [\ref{CS3}] generalize straightforwardly
to meromorphic $T$-actions.
It is easy to see that the following statements are equivalent:\\
1. $(F,\prec)$ is a partially ordered set;\\
2. There is no quasicycle in $(F,\prec)$;\\
3. There is a strictly decreasing function on $(F,\prec)$, i.e., a function
$f\colon F\to\RE$ satisfying $f(\al)>f(\beta)$ if $\al\prec\beta$ and 
$\al\ne\beta$.\\
Consequently, $(F,\prec)$ is a partially ordered set if one of the following
is true:\\
1. $X$ is K\"ahler;\\
2. $\nu^C_\al>\nu^C_\beta$ if $\al\prec\beta$ and $\al\ne\beta$;\\
3. The \bb decomposition is a stratification.\\
In each of the above cases, the moment map (projected along some direction
in $C$), $\nu^C_\cdot$, and $\dim_\co X^C_\cdot$, respectively,
provides a strictly decreasing function on $(F,\prec)$.}
\end{rmk}

\begin{ex}\label{JW}
{\em Jurkiewicz [\ref{J}] constructed a smooth compact toric $3$-manifold 
with a meromorphic $T^3$-action that has $22$ isolated fixed points.
Choosing an appropriate action chamber, there is a quasicycle of length $6$
[\ref{J}].
Therefore $(F,\prec)$ is not a partially ordered set.
In [\ref{Wu}, \S\,4], it is shown that there exists a $T^3$-equivariant
holomorphic line bundle such that the holomorphic Morse inequalities fail.
This shows that the holomorphic Morse inequalities are not valid on an 
arbitrary complex manifold [\ref{Wu}], though the fixed-point theorems in
[\ref{AB}, \ref{AS}] requires no further assumptions.
In section 3, we construct the analog of the instanton complex in 
holomorphic Morse theory when $(F,\prec)$ is a partially ordered set.
The existence of such a construction implies the holomorphic Morse
inequalities.
Moreover, the partial order condition is weaker than the \ka condition.}
\end{ex}

\begin{defn}
{\em Suppose $X$ has a $C$-meromorphic $T$-action.
The \bb decomposition with respect to $C$ is {\em filterable} if there is
a descending sequence of $T$-invariant subvarieties
	\be\label{closed}
X=Z_0\supset Z_1\supset\cdots\supset Z_m\supset Z_{m+1}=\emptyset
	\ee
such that for all $0\le p\le m$, $Z_p-Z_{p+1}=\bigcup_{\al\in F_p}X^C_\al$
for a subset $F_p\subset F$ such that neither $\al\prec\beta$ nor
$\beta\prec\al$ if $\al\ne\beta\in F_p$.}
\end{defn}

Notice that we allow $Z_p-Z_{p+1}$ to be a union of cells labeled by
elements in $F$ unrelated by $\prec$.
In [\ref{BB2}, Definition~2], $Z_p-Z_{p+1}$ is required to be a single cell.
Since $\overline{X^C_\al}\bigcap X^C_\beta\ne\emptyset$ implies
$\al\prec\beta$ [\ref{CS3}, Lemma~1], the two notions are equivalent.
Notice that the function $\al\mapsto p(\al)$ where $\al\in F_{p(\al)}$ is
strictly increasing on $(F,\prec)$.

Alternatively, (\ref{closed}) can be written as
	\be\label{open}
X=V_0\supset V_1\supset\cdots\supset V_m\supset V_{m+1}=\emptyset,
	\ee
where $V_p=X-Z_{m+1-p}$ ($0\le q\le m+1$) are open sets in $X$ such that
$V_p-V_{p+1}=Z_{m-p}-Z_{m-p+1}=\bigcup_{\al\in F_{m-p}}X^C_\al$ 
for $0\le p\le m$.

We return to the simple example $X=\CP^1$ with two fixed points $0$, $\infty$
under the meromorphic $\CM$-action.
Let $Z_0=V_0=X$, $Z_1=\{\infty\}$, $V_1=X-\{\infty\}$, $Z_2=V_2=\emptyset$,
then $Z_0-Z_1=V_1-V_2=X^+_0$, $Z_1-Z_2=V_0-V_1=X^+_\infty$.
So the plus decomposition of $X=\CP^1$ is filterable.
For the same reason, so is the minus decomposition.

\begin{prop}\label{CELL} {\em ([\ref{CS3}])}
Consider a meromorphic $T$-action on $X$.
Then the \bb decomposition {\em (\ref{bbc})} is filterable if and only if
$(F,\prec)$ is a partially ordered set.
If so, then\\
1. The projection $\pi^C\colon X^C_\al\to X^T_\al$ is a $T$-equivariant
holomorphic fibration and the fiber $(\pi^C)^{-1}(x)$ over any $x\in X^T_\al$
is $T$-equivariantly isomorphic to $(N^C_\al)_x$;\\
2. There is a $T$-equivariant isomorphism 
$TX^C_\al|_{X^T_\al}\cong N^C_\al\oplus TX^T_\al$
of holomorphic vector bundles over $X^T_\al$;\\
3. The closure $\overline{X^C_\al}$ in $X$ is analytic and contains 
$X^C_\al$ as a Zariski open set. 
Consequently, $X^C_\al$ is locally closed in $X$.
\end{prop}

\proof{If $T=\CM$, the necessary and sufficient condition for (\ref{bbc})
to be filterable was proved in [\ref{CS3}].
Properties 1 and 2 follow from the arguments of [\ref{CS1}].
Property 3 follows from the arguments in [\ref{CS2}, \S\,IIb], where the \ka
assumption was not made.
The generalization to a higher rank torus $T$ is straightforward.}

\vspace{1ex}
The three properties of Proposition~\ref{CELL} were shown to be valid 
when $X$ is a \ka manifold [\ref{CS1}, \ref{CS2}, \ref{Fu}, \ref{Kor}] 
or a complete normal algebraic variety [\ref{BB1}, \ref{BB2}, \ref{Kon}],
prior to the work of [\ref{CS3}].
Without any of these assumptions, one or more of the properties
in Proposition~\ref{CELL} could fail [\ref{So}].

Example \ref{JW} was originally constructed to provide
a non-filterable \bb decomposition [\ref{J}]. 

\begin{rmk}
{\em The restriction of $\pi^{-C}$ to $X^C_\al-X^T_\al$ may be discontinuous
and the image $\pi^{-C}(X^C_\al-X^T_\al)$ may fall into more than one
connected components of $X^T$.
For example, let $X=\CP^1\times\CP^1$ with the diagonal $\CM$-action.
Then $X^T=\{0,\infty\}\times\{0,\infty\}$ and $X^+_{(0,0)}=\CO\times\CO$.
We have $\pi^-(\zero\times(\CO-\zero))=(0,\infty)$,
$\pi^-((\CO-\zero)\times\zero)=(\infty,0)$, and
$\pi^-((\CO-\zero)\times(\CO-\zero))=(\infty,\infty)$.
The reason is that the holomorphic embedding $X^C_\al\to X$ extend 
only meromorphically at infinity [\ref{CS2}, Lemma~2], where it can 
be discontinuous.}
\end{rmk}

Notice that despite of part~2 of Proposition~\ref{CELL},
a tubular neighborhood of $X^T_\al$ in $X^C_\al$ can not be identified 
holomorphically with that in $N^C_\al$ in general [\ref{CS1}].
There is an infinite series of obstruction to this [\ref{Gr}, \ref{EL}].
However, an identification is possible locally on $X^T_\al$.
Consider a holomorphic vector bundle $E$ over $X$ on which the $T$-action 
lifts holomorphically.
For future applications, we also put $E$ into a standard local form.

\begin{lemma}\label{LOCAL}
For any $x\in X^T_\al$, there is a neighborhood $U_x$ of $x$ in $X^T_\al$,
a $T_\re$-invariant open set $W^C_x$ in $N_\al$ containing $N^C_\al|_{U_x}$
as a closed subset,
and a $T$-equivariant holomorphic embedding $\psi_x\colon W^C_x\to X$
such that $\psi_x(N^C_\al|_{U_x})=(\pi^C)^{-1}(U_x)\subset X^C_\al$.
Moreover, $\psi_x$ can be lifted to a $T$-equivariant isomorphism
$\tilde{\psi}_x\colon W^C_x\times E_x\to E|_{\psi_x(W^C_x)}$
of holomorphic vector bundles.
\end{lemma}

\proof{As in the proof of [\ref{CS1}, Proposition~I], there is a neighborhood
$U_x$ of $x$ in $X^T_\al$, a $T_\re$-invariant open set $W_x$ in $N_\al$ 
containing $U_x$, and a $T$-equivariant holomorphic embedding 
$\psi_x\colon W_x\to X$ such that $\psi_x(N^C_\al\cap W_x)\subset X^C_\al$.
Pick any $v\in\ell\cap C$.
Let $W^C_x=\bigcup_{t\ge0}j_v(\e{t})\,W_x$.
$W^C_x$ is a $T_\re$-invariant open set in $N_\al$.
Moreover, for any $y\in N^C_\al|_{U_x}$, we have $\pi^C(y)\in U_x$, hence 
there exists $t\ge0$ such that $j_v(\e{-t})y\in W_x$, i.e., $y\in W^C_x$.
So $W^C_x$ contains $N^C_\al|_{U_x}$.
We extend $\psi_x$ from $W_x$ to $W^C_x$ by 
$\psi_x(j_v(\e{t})y)=j_v(\e{t})\psi_x(y)$ for $y\in W_x$ and $t\ge0$.
Clearly, the extension is well-defined, $T$-equivariant and holomorphic.
Next, there is a holomorphic isomorphism
$\tilde{\psi}_x\colon W_x\times E_x\to E|_{W_x}$ of vector bundles,
perhaps on a smaller neighborhood $W_x$.
By [\ref{CS1}, Lemma~I], $\tilde{\psi}_x$ can be made $T_\re$-equivariant
(hence $T$-equivariant).
We extend $\tilde{\psi}_x$ to $W^C_x\times E_x$ by
$\tilde{\psi}_x(j_v(\e{t})y,\xi)=j_v(\e{t})\tilde{\psi}_x(y,j_v(\e{-t})\xi)$ 
for $y\in W_x$, $\xi\in E_x$ and $t\ge0$.
The extension is again well-defined and is a $T$-equivariant holomorphic
isomorphism of vector bundles.}

\subsect{A non-compact setting}

In this subsection, we consider a class of non-compact complex manifolds
with holomorphic torus actions.
We hope that this class is broad enough to include many interesting examples.

Let $X$ be a (possibly non-compact) complex manifold with a holomorphic
action of the torus $T$.

\begin{defn}\label{CMERO} 
{\em Let $C$ be an action chamber.
The $T$-action on $X$ is {\em $C$-meromorphic} if\\
1. for any $x\in X$, $v\in\ell\cap C$, the limit $\pi^v(x)$ exists;\\
2. there is a compact complex orbifold $\tilde{X}$ with a meromorphic 
$T$-action and a $T$-equivariant holomorphic embedding of $X$ onto 
a Zariski open set of $\tilde{X}$.}
\end{defn}

The simplest example is $X=\CO$ with the standard multiplication by $\CM$.
The action is plus-meromorphic and $X$ has a compactification 
$\tilde{X}=\CP^1$.
The plus-decomposition is $X=X^+_0$ (a single $2$-cell) and there is no
minus-decomposition.

\begin{rmk} 
{\em Consider a $C$-meromorphic $T$-action on $X$.
We identify $X$ with its image in $\tilde{X}$.\\
1. By Proposition~\ref{LIMIT}, which applies to the non-compact setting here,
the limit $\pi^v(x)$ ($x\in X$) does not depend on the choice of 
$v\in\ell\cap C$ and is therefore denoted by $\pi^C(x)$.
Moreover $\pi^C(x)\in X^T$.
Because $X$ is embedded into a compact space $\tilde{X}$, the set $F$ 
of connected components of $X^T$ is finite and each component $X^T_\al$ 
($\al\in F$) is compact.
We have the \bb decomposition (\ref{bbc}) with respect to $C$.
The action chamber $C$ of $X$ may be divided into several action chambers
of $\tilde{X}$; let $\tilde{C}$ be one of such.
Then we have $X^C_\al=\tilde{X}^{\tilde{C}}_\al$ for any $\al\in F$.
For $x\in X$, the limit $\pi^{-\tilde{C}}(x)$ exists in $\tilde{X}$
but may fall into $\tilde{X}-X$.
Therefore $\pi^{-C}(x)$ is in general not defined in $X$.\\
2. As in the compact situation, there is a relation $\prec$ on $F$.
If $(F,\prec)$ is a partially ordered set, then the properties of 
Proposition~\ref{CELL} for $X$ are satisfied.
In particular, $\overline{X^C_\al}=\overline{\tilde{X}^{\tilde{C}}_\al}\cap X$
is a closed subvariety in $X$ that contains $X^C_\al=\tilde{X}^{\tilde{C}}_\al$
as a Zariski open set.
Furthermore, the \bb decomposition of $X$ with respect to $C$ is filterable
and we have filtrations of $X$ by closed subsets (\ref{closed}) and by open
subsets (\ref{open}).\\
3. If $E$ is a holomorphic vector bundle over $X$ on which the $T$-action
lifts holomorphically, Lemma~\ref{LOCAL} also holds.}
\end{rmk}

\begin{assump}\label{WA}{\em
There exists an action chamber $C$ such that the $T$-action on $X$ is
$C$-meromorphic and the set $(F,\prec)$ is partially ordered.}
\end{assump}

In section 3, we will establish holomorphic Morse theory on a (possibly
non-compact) complex manifolds satisfying Assumption~\ref{WA}.
An immediate way of obtaining such non-compact manifolds comes from
Definition~\ref{CMERO}.
We start with a compact complex manifold $\tilde{X}$ with a meromorphic
$T$-action.
Suppose the \bb decomposition of $\tilde{X}$ with respect to an action 
chamber $\tilde{C}$ is filterable and is filtered by the closed sets 
	\be
\tilde{X}=\tilde{Z}_0\supset\tilde{Z}_1\supset\cdots\supset
\tilde{Z}_{\tilde{m}}\supset\tilde{Z}_{\tilde{m}+1}=\emptyset.
	\ee
Pick any $m$ such that $0\le m\le\tilde{m}-1$ and 
let $X=\tilde{X}-\tilde{Z}_{m+1}$.
$T$ acts holomorphically on $X$.
Let $C$ be the action chamber that contains $\tilde{C}$.
Then the $T$-action on $X$ is $C$-meromorphic.
Moreover the \bb decomposition of $X$ with respect to $C$ has a filtration
(\ref{closed}) by closed subsets $Z_p=\tilde{Z}_p-\tilde{Z}_{m+1}$ 
($0\le p\le m+1$) of $X$.
The simple example $X=\CO$ falls into this category, with $\tilde{X}=\CP^1$.

More interestingly, the non-compact setting here is a complex analog of 
the symplectic setting considered in [\ref{P}, \ref{PW}], which we now recall.
Let $(X,\om)$ be a (possibly non-compact) symplectic manifold with
a Hamiltonian action of the compact torus $T_\re$, with a moment map
$\mu\colon X\to\gt^*_\re$.
The fixed-point set $X^T$ of the torus $T_\re$ is a symplectic submanifold 
of $X$.
Let $F$ be the set of connected components of $X^T$.

\begin{assump}\label{PWA}
{\em ([\ref{PW}, Assumption~1.3])
There is $v\in\gt_\re$ such that $\bra\mu,v\ket\colon X\to\RE$ is proper 
and not surjective and $F$ is a (non-empty) finite set.}
\end{assump}

If in addition $(X,\om)$ is \ka and the $T_\re$-action preserves the complex
structure on $X$, then there is a holomorphic $T$-action on $X$.

\begin{prop}
Let $(X,\om)$ be a \ka manifold with a holomorphic $T$-action.
Suppose that the $T_\re$-action is Hamiltonian.
Then Assumption~\ref{PWA} implies Assumption~\ref{WA}.
\end{prop}

\proof{By [\ref{PW}, Proposition~1.6], there is an action chamber $C$ 
such that for any $v\in C$, the function $\bra\mu,v\ket$ on $X$ is proper
and bounded from above.
Therefore if $v\in\ell\cap C$, the limit $\pi^v(x)$ exists for any $x\in X$.
Pick any $v\in\ell\cap C$.
Since $F$ is finite, there is $a\in\RE$ such that $\bra\mu(X^T),v\ket>a$.
We construct a symplectic cut $X_{\ge a}$ [\ref{Le}].
Let $\CM$ act on $X\times\CO$ by $u\colon(x,z)\mapsto(j_v(u)x,uz)$.
The action of $S^1\subset\CM$ on $X\times\CO$ is Hamiltonian with a moment map
$\tilde\mu(x,z)=\bra\mu(x),v\ket-a-\hf|z|^2$.
$\tilde\mu$ is a proper function on $X\times\CO$ and $0$ is a regular value.
The symplectic quotient $X_{\ge a}=\tilde\mu^{-1}(0)/S^1$ is a compact
symplectic orbifold with a Hamiltonian $T_\re$-action.
Since $X$ is K\"ahler, $X_{\ge a}=(X\times\CO)^{\rm s}/\CM$ holomorphically
and is also K\"ahler [\ref{GS}].
Here $(X\times\CO)^{\rm s}=\set{(x,z)\in X\times\CO}
{\CM(x,z)\cap\tilde\mu^{-1}(0)\ne\emptyset}$ is the stable subset of
$X\times\CO$.
We want to construct a $T$-equivariant holomorphic embedding $X\to X_{\ge a}$.
Clearly, $\tilde\mu(u(x,1))=\bra\mu(j_v(u)x),v\ket-a-\hf|u|^2$,
where $u\in\CM$ and $x\in X$.
For any $x\in X$, since $\pi^v(x)\in X^T$, 
$\lim_{u\to0}\tilde\mu(u(x,1))=\bra\mu(\pi^v(x)),v\ket-a>0$.
On the other hand, since $\bra\mu,v\ket$ is bounded from above,
$\lim_{u\to\infty}\tilde\mu(u(x,1))=-\infty$.
Therefore there is $u\in\CM$ such that $\tilde\mu(u(x,1))=0$.
Hence $X\times\{1\}\subset(X\times\CO)^{\rm s}$.
The composition $X\to X\times\{1\}\subset(X\times\CO)^{\rm s}\to X_{\ge a}$
of the inclusion and the quotient is a $T$-equivariant holomorphic embedding.
The image is $X_{>a}=(X\times(\CO-\zero))/\CM$.
Since $X_{\ge a}-X_{>a}=(X\times\zero)^{\rm s}/\CM=\mu^{-1}(a)/S^1$ is 
a complex subvariety of $X_{\ge a}$, $X$ is embedded as a Zariski open set.
Using the moment map $\bra\mu,v\ket$, it is easy to show that $(F,\prec)$
is a partially ordered set.}

\sect{Equivariant spectral sequences in holomorphic Morse theory}

We consider a holomorphic $T$-action on a (possibly non-compact) 
complex manifold $X$.
$F$ is the set of connected components of the fixed-point set $X^T$.
Throughout this section, we make Assumption~\ref{WA}.
Then the \bb decomposition is filterable, with descending sequences of
closed sets (\ref{closed}) and open sets (\ref{open}) in $X$.
Let $E$ be a holomorphic vector bundle over $X$ on which the $T$-action
lifts holomorphically.
We want to determine the Dolbeault cohomology groups $H^*\cpct(X,\OO(E))$ 
(with compact support) and $H^*(X,\OO(E))$ as representations of $T$.

\begin{defn}
{\em A (cohomological) spectral sequence $\{E\pqr, \dr\pqr\}$ is}
$T$-equivariant {\em if the spaces $E\pqr$ are representations of $T$ and the
coboundary maps $\dr\pqr\colon E\pqr\to E^{p+r,q-r+1}_r$ are $T$-equivariant.
The spectral sequence} converges $T$-equivariantly {\em to the representations
$H^*$ if the spaces $E\pq_\infty$ are the graded components of $H^*$ 
as representations of $T$.}
\end{defn}

\subsect{Spectral sequence for cohomologies with compact support}

In this subsection, we construct a spectral sequence converging to 
the Dolbeault cohomology groups $H^*\cpct(X,\OO(E))$ with compact support.

Recall that if $A\subset X$ is a locally closed subset, then for any sheaf
$\FF$ on $X$, there is a unique sheaf on $X$, denoted by $\FF_A$, such that
the restrictions $\FF_A|_A=\FF|_A$ and $\FF|_{X-A}=0$.
Moreover, $\FF_A$ exists for any sheaf $\FF$ only if $A$ is locally closed
[\ref{G}, Th\'eor\`eme II.2.9.1].
Let $0\to\FF\to\CF$ be the canonical resolution of $\FF$ [\ref{G}, \S\,II.4.3].
It is easy to see that $0\to\FF_A\to\CF_A$ is a flabby resolution of $\FF_A$.
Finally, if $A$ is an open subset, then $\FF_A$ is a subsheaf of $\FF$.

For simplicity, we denote the sheaf $\OO(E)$ by $\FF$ from now on.
If $A$ is a $T$-invariant locally closed subset of $X$, then $T$ acts on 
the sheaf $\FF_A$ and hence on the cohomology groups $H^*\cpct(X,\FF_A)$.

\begin{lemma}\label{FILT}
Under Assumption~\ref{WA}, there is a $T$-equivariant spectral sequence with
	\be
E\pq_1=H^{p+q}\cpct(X,\FF_{V_p-V_{p+1}})
	\ee
that converges $T$-equivariantly to $H^*\cpct(X,\FF)$.
\end{lemma}

\proof{From (\ref{open}), we have a filtration of $\FF$ by subsheaves
	\be\label{fsheaf}
\FF=\FF_{V_0}\supset\FF_{V_1}\supset\cdots
\supset\FF_{V_m}\supset\FF_{V_{m+1}}=0
	\ee
and hence a filtration of the cochain complex $\Gamc(\CF)$ by
	\be
\Gamc(\CF)=\Gamc(\CF_{V_0})\supset\Gamc(\CF_{V_1})\supset\cdots
\supset\Gamc(\CF_{V_m})\supset\Gamc(\CF_{V_{m+1}})=0.
	\ee
This induces a spectral sequence that converges to 
$H^*(\Gam(\CF))=H^*(X,\FF)$, with
	\be
E\pq_0=\Gamc(\CC^{p+q}(\FF)_{V_p})/\Gamc(\CC^{p+q}(\FF)_{V_{p+1}})
=\Gamc(\CC^{p+q}(\FF)_{V_p-V_{p+1}}).
	\ee
Since the maps $d\pq_0\colon E\pq_0\to E^{p,q+1}_0$ are induced 
by the resolution, we get
	\be
E\pq_1=H^{p+q}(\Gamc(\CC^*(\FF)_{V_p-V_{p+1}}))=
H^{p+q}\cpct(X,\FF_{V_p-V_{p+1}}).
	\ee
All the steps are $T$-equivariant.}

\begin{lemma}\label{SPLIT}
	\be\label{split}
H^*\cpct(X,\FF_{V_p-V_{p+1}})=
\bigoplus_{\al\in F_{m-p}}H^*\cpct(X^C_\al,\FF|_{X^C_\al})
	\ee
as representations of $T$.
\end{lemma}

\proof{Since $\overline{X^C_\al}\cap X^C_\beta=\emptyset$ 
for any $\al\ne\beta\in F_{m-p}$, we have 
$\FF_{V_p-V_{p+1}}=\bigoplus_{\al\in F_{m-p}}\FF_{X^C_\al}$ and hence
	\be\label{step1}
H^*\cpct(X,\FF_{V_p-V_{p+1}})=
\bigoplus_{\al\in F_{m-p}}H^*\cpct(X,\FF_{X^C_\al}).
	\ee
The support of $\FF_{X^C_\al}$ is contained in the closed subvariety 
$\overline{X^C_\al}$.
Therefore we have [\ref{G}, Th\'eor\`em II.4.9.1]
	\be\label{step2}
H^*\cpct(X,\FF_{X^C_\al})=H^*\cpct(\overline{X^C_\al},\FF_{X^C_\al}).
	\ee
Since $\overline{X^C_\al}-X^C_\al$ is a closed subset in $\overline{X^C_\al}$
and $\FF_{X^C_\al}|_{\overline{X^C_\al}-X^C_\al}=0$, we deduce from
[\ref{G}, Th\'eor\`em~II.4.10.1] that
	\be\label{step3}
H^*\cpct(X^C_\al,\FF|_{X^C_\al})=H^*\cpct(\overline{X^C_\al},\FF_{X^C_\al}).
	\ee
The result follows from (\ref{step1}), (\ref{step2}) and (\ref{step3}).}

\vspace{1ex}
Recall that $\pi^C\colon X^C_\al\to X^T_\al$ is a holomorphic fibration
with fiber $\CO^{\nu^C_\al}$.
The sheaf $\FF|_{X^C_\al}$ is on the total space $X^C_\al$.
To calculate the right hand side of (\ref{split}), we need another
spectral sequence.

We consider a general fibration $\pi\colon Y\to B$ over a compact base $B$
with possibly non-compact fibers.
For the time being, let $\FF$ be an arbitrary sheaf on the total space $Y$.
The cohomology groups with compact support are 
$H^q\cpct(Y,\FF)=H^q(\Gam_\PHI(Y,\CF))$ ($q\ge0$), where $\PHI$ is
a family of supports that consists of the compact subsets of $Y$.
Let $\A$, $\LL^*$ be the sheaves on $B$ defined by the presheaves
$\A(U)=\Gam_{\PHI\cap\pi^{-1}(U)}(\pi^{-1}(U),\FF)$,
$\LL^*(U)=\Gam_{\PHI\cap\pi^{-1}(U)}(\pi^{-1}(U),\CF)$, respectively,
where $U$ is any open subset of $B$.
Then $0\to\A\to\LL^*$ is a differential sheaf in the sence of
[\ref{G}, \S\,II.4.1].
Let $\HH^q\cpct(Y,\FF)$ ($q\ge0$) be the sheaves on $B$ 
defined by the presheaves $\HH^q\cpct(Y,\FF)(U)=H^q(\LL^*(U))$, for any open
subset $U\subset B$.

\begin{lemma}\label{PROJ}
1. At $b\in B$, the stalk of $\HH^q\cpct(Y,\FF)$ for any $q\ge0$ is
	\be
\HH^q\cpct(Y,\FF)_b\cong H^q\cpct(Y_b,\FF|_{Y_b}).
	\ee
2. There is a spectral sequence with
	\be\label{e2}
E\pq_2=H^p(B,\HH^q\cpct(Y,\FF))
	\ee
that converges to $H^*\cpct(Y,\FF)$.
\end{lemma}

\proof{1. This the analog of [\ref{G}, Remarque II.4.17.1] for cohomologies
with compact support.
First, $\HH^q\cpct(Y,\FF)_b=\indlim_{U\ni b}\HH^q\cpct(Y,\FF)(U)=
\indlim_{U\ni b}H^q_{\PHI\cap\pi^{-1}(U)}(\pi^{-1}(U),\FF)$.
By [\ref{G}, Th\'eor\`em II.3.3.1], any section $s\in\Gam(Y_b,\CF|_{Y_b})$
can be extended to a neighborhood of $Y_b$ in $Y$.
If $\supp s\in\PHI\cap Y_b$, then the neighborhood can be chosen as 
$\pi^{-1}(U)$ for some open set $U\subset B$.
Therefore $\indlim_{U\ni b}H^q_{\PHI\cap\pi^{-1}(U)}(\pi^{-1}(U),\FF)=
\indlim_{V\supset Y_b}H^q_{\PHI\cap V}(V,\FF)$.
Following the proof of [\ref{G}, Th\'eor\`em II.4.11.1], we get
$\indlim_{V\supset Y_b}H^q_{\PHI\cap V}(V,\FF)=H^q\cpct(Y_b,\FF|_{Y_b})$.\\
2. It is clear that $\LL^*$ are flabby sheaves.
By [\ref{G}, Th\'eor\`em II.4.6.1], associated to the differential sheaf
$0\to\A\to\LL^*$ there is a spectral sequence with (\ref{e2}) that converges
to $H^*(\Gam(\LL^*))$.
By the definition of $\LL^*$, $\Gam(B,\LL^*)=\Gam\cpct(Y,\CF)$.
The result follows.}

\begin{lemma}\label{MODEL}
Let $R$ be a representation of $T$ with $\dim_\co R=n$ and let $A$, $A^\perp$
be $T$-invariant subspaces of $R$ such that $\dim_\co A=\nu$ and 
$R=A\oplus A^\perp$.
Let $E_0$ be a representation of $T$ and $E$, the trivial holomorphic vector
bundle over $V$ with fiber $E_0$.
Then, as representations of $T$,
	\be
H^q\cpct(A,\OO(E)|_A)=H^q\cpct(R,\OO(E)_A)=
\two{S(A^{\perp*})\times S(A)\otimes\medwedge^\nu(A)\otimes E_0}{q=\nu}
{0}{q\ne\nu.}
	\ee
\end{lemma}

\proof{It suffices to prove the case when $E_0=\CO$ is a trivial 
representation.
If $A=\zero$, then
	\be
H^q\cpct(R,\OO_\zero)=H^q(R,\OO_\zero)=\two{S(R^*)}{q=0}{0}{q\ne0.}
	\ee
If $A=R$, then (see [\ref{La}] for an analytic version)
	\be
H^q\cpct(R,\OO)=\two{S(R)\otimes\medwedge^n(R)}{q=n}{0}{q\ne n.}
	\ee
The general case is a consequence of the K\"unneth formula.}

\vspace{1ex}
We now return to the situation of $\FF=\OO(E)$.

\begin{lemma}\label{EACH}
	\be\label{each}
H^q\cpct(X^C_\al,\FF|_{X^C_\al})=H^{q-\nu^C_\al}(X^T_\al,\OO(S((N^{-C}_\al)^*)
\otimes S(N^C_\al)\otimes\medwedge^{\nu^C_\al}(N^C_\al)\otimes E|_{X^T_\al}))
	\ee
as representations of $T$.
\end{lemma}

\proof{Consider the holomorphic fibration $\pi^C\colon X^C_\al\to X^T_\al$
with fiber $\CO^{\nu^C_\al}$ and the sheaf $\FF|_{X^C_\al}$ on $X^C_\al$.
For any $x\in X^T_\al$, we want to find the stalk $\HH^q\cpct(N^C,\FF)_x$,
which depends only on an open neighborhood of $(\pi^C)^{-1}(x)\subset X^C_\al$
in $X$.
By Lemma~\ref{LOCAL}, we can replace $X^C_\al\subset X$ by 
$N^C_\al|_{U_x}\subset N_\al|_{U_x}$ and $E$ by a trivial vector bundle
with fiber $E_x$.
Moreover there is a $T$-equivariant isomorphism
$(N_\al,N^C_\al)|_{U_x}\cong U_x\times(N_x,N^C_x)$.
By Lemma~\ref{PROJ}.1 and Lemma~\ref{MODEL},
	\bea
\HH^q\cpct(N^C,\FF)_x\eq H^q(N^C_x,\OO(W^C_x,E_x)|_{N^C_x})		\nno
\eq\two{\OO(S((N^{-C}_\al)^*)\otimes S(N^C_\al)\otimes
\medwedge^{\nu^C_\al}(N^C_\al)\otimes E|_{X^T_\al})_x}
{q=\nu^C_\al}{0}{q\ne\nu^C_\al.}
	\eea
So the spectral sequence of Lemma~\ref{PROJ}.2 degenerates at $E_2$
and the result follows.}

Though the bundle $S((N^{-C}_\al)^*)\otimes S(N^C_\al)\otimes
\medwedge^{\nu^C_\al}(N^C_\al)\otimes E|_{X^T_\al}$ over $X^T_\al$
is infinite dimensional, its sub-bundle of any given weight is of finite rank.
Therefore each weight has a finite multiplicity in the cohomology groups
(\ref{each}), and their formal characters in $\ZZ[\ell^*]$ exist.

\begin{thm}\label{MAIN}
Let $X$ be a complex manifold with a holomorphic $T$-action satisfying
Assumption~\ref{WA}.
Let $E$ be a holomorphic vector bundle over $X$ on which the $T$-action lifts
holomorphically.
Then\\
1. there is a $T$-equivariant spectral sequence converging $T$-equivariantly
to $H^*\cpct(X,\OO(E))$ with
	\be\label{e1}
E\pq_1=
\bigoplus_{\al\in F_{m-p}}H^{p+q-\nu^C_\al}(X^T_\al, \OO(S((N^{-C}_\al)^*)
\otimes S(N^C_\al)\otimes\medwedge^{\nu^C_\al}(N^C_\al)\otimes E|_{X^T_\al}));
	\ee
2. there is a character valued polynomial $Q^C\cpct(t)\ge0$ such that
	\bea\label{morse}
\vc\sum_{\al\in F}t^{\nu^C_\al}\sum_{q=0}^{n_\al}t^q
\ch H^q(X^T_\al,\OO(S((N^{-C}_\al)^*)\otimes S(N^C_\al)\otimes
\medwedge^{\nu^C_\al}(N^C_\al)\otimes E|_{X^T_\al}))		\nno
\eq\sum_{q=0}^nt^q\ch H^q\cpct(X,\OO(E))+(1+t)Q^C\cpct(t);
	\eea
3. 
	\be\label{ind}
\sum_{q=0}^n(-1)^q\ch H^q\cpct(X,\OO(E))=
\sum_{\al\in F}(-1)^{\nu^C_\al}\int_{X^T_\al}
{\mathrm{ch}}_T\left(\frac{E|_{X^T_\al}\otimes\det(N^C_\al)}
{\det(1-(N^{-C}_\al)^*)\otimes\det(1-N^C_\al)}\right){\mathrm{td}}(X^T_\al),
	\ee
where $\mathrm{ch}_T$ and $\mathrm{td}$ stand for the equivariant Chern 
character and the Todd class, respectively.
\end{thm}

\proof{1. The result follows from Lemma~\ref{FILT}, Lemma~\ref{SPLIT} 
and Lemma~\ref{EACH}.\\
2. Since $E\pq_{r+1}$ is the cohomology of $(E\pqr,\dr\pqr)$, we have
	\be\label{r+1}
\sum_{p,q}t^{p+q}\ch E\pqr=\sum_{p,q}t^{p+q}\ch E\pq_{r+1}+(1+t)Q_r(t)
	\ee
for a character valued polynomial $Q_r(t)\ge0$.
Using (\ref{r+1}) recursively, we get (\ref{morse}) with 
$Q^C\cpct(t)=\sum_{r\ge1}Q_r(t)\ge0$.\\
3. By setting $t=-1$ in (\ref{morse}) and using 
	\bea
\vc\sum_{q=0}^{n_\al}(-1)^q\ch H^q(X^T_\al,\OO(S((N^{-C}_\al)^*)\otimes 
S(N^C_\al)\otimes\medwedge^{\nu^C_\al}(N^C_\al)\otimes E|_{X^T_\al}))	\nno
\eq\int_{X^T_\al}{\mathrm{ch}}_T\left(\frac{E|_{X^T_\al}\otimes\det(N^C_\al)}
{\det(1-(N^{-C}_\al)^*)\otimes\det(1-N^C_\al)}\right){\mathrm{td}}(X^T_\al),
	\eea
we obtain (\ref{ind}).
See [\ref{WZ}, Remark 2.3.2].}

\begin{cor}\label{DISCR}
If in addition $X^T$ is discrete {\em ({\em and is identified with $F$})},
then\\
1. there is a $T$-equivariant spectral sequence converging $T$-equivariantly to
$H^*\cpct(X,\OO(E))$ with
	\be\label{e1-discr}
E\pq_1=\bigoplus_{x\in F,\,\nu^C_x=p+q}
S((N^{-C}_x)^*)\otimes S(N^C_x)\otimes\medwedge^{\nu^C_x}(N^C_x)\otimes E_x.
	\ee
2. there is a character valued polynomial $Q^C\cpct(t)\ge0$ such that
	\be\label{morse-discr}
\sum_{x\in F}t^{\nu^C_x}\ch E_x
\prod_{\lam_{x,k}\in C^*}\frac{\e{\lam_{x,k}}}{1-\e{\lam_{x,k}}}
\prod_{\lam_{x,k}\in-C^*}\inv{1-\e{-\lam_{x,k}}}
=\sum_{q=0}^nt^q\ch H^q\cpct(X,\OO(E))+(1+t)Q^C\cpct(t);
	\ee
3. 
	\be\label{ind-discr}
\sum_{q=0}^n(-1)^q\ch H^q\cpct(X,\OO(E))=\sum_{x\in F}(-1)^{\nu^C_\al}
\ch E_x\prod_{\lam_{x,k}\in C^*}\frac{\e{\lam_{x,k}}}{1-\e{\lam_{x,k}}}
\prod_{\lam_{x,k}\in-C^*}\inv{1-\e{-\lam_{x,k}}}.
	\ee
\end{cor}

\begin{rmk}\label{REMARK}{\em
1. If $X$ is compact, then $H^*\cpct(X,\OO(E))=H^*(X,\OO(E))$, and the
right-hand sides of (\ref{ind}) and (\ref{ind-discr}) are often written as
	\be
\sum_{\al\in F}\int_{X^T_\al}{\mathrm{ch}}_T\left(
\frac{E|_{X^T_\al}}{\det(1-N_\al^*)}\right){\mathrm{td}}(X^T_\al)
\quad\quad\mbox{and}\quad\quad
\sum_{x\in F}\frac{\ch E_x}{\prod_{k=1}^n(1-\e{-\lam_{x,k}})},
	\ee
respectively.
In this case, parts~3 of Theorem~\ref{MAIN} and Corollary~\ref{DISCR}
are the fixed-point theorems of [\ref{AS}, \ref{AB}],
which do not require $\prec$ to be a partial ordering.
Here $X$ can be non-compact.
We obtain a fixed-point theorem for Dolbeault cohomology groups
with compact support under the partial order condition.
When $X$ is compact and K\"ahler, parts~2 are the results of 
[\ref{W2}, \ref{MW}, \ref{Wu}, \ref{WZ}].
Parts~1 strengthen these results under a weaker condition, namely, 
Assumption~\ref{WA}.
In particular, all the weights of $T$ in $H^*\cpct(X,\OO(E))$ are of
finite multiplicity.
It would be interesting to have an independent analytic proof of the results
in parts~2 when $X$ is a non-compact K\"ahler manifold satisfying
Assumption~\ref{PWA}.
They are the discrete versions of [\ref{PW}, Theorem~3.2].\\
2. The coboundary maps $\{\dr\pqr\}$ in the spectral sequence in
Theorem~\ref{MAIN} or Corollary~\ref{DISCR} are the holomorphic counterparts
of the instanton tunneling operators in [\ref{W1}].
Through this spectral sequence, the cohomology groups $H^*\cpct(X,\OO(E))$ are
completely determined by the combinatorial data of the $T$-action on $X$.
However unlike the real case, the spectral sequence of holomorphic 
Morse theory does not always degenerate at $E_2$.
A sufficient condition for degeneracy at $E_2$ is
	\be\label{dege2}
E\pq_1=0\quad\mbox{for all}\quad q\ne0.
	\ee
If so, then the spectral sequence reduces to a cochain complex
$\{E^{*0}_1,\dr^{*0}_1\}$, whose cohomology is $E^{*0}_2=H^*(X,\OO(E))$.
This would be exactly like the Thom-Smale-Witten complex [\ref{W1}].
For example, if $X^T=F$ is discrete, $m=n$ in (\ref{open}), and 
$F_p=\set{x\in F}{\nu^{-C}_x=p}$ for $0\le p\le n$,
then (\ref{dege2}) is satisfied.}
\end{rmk}

\subsect{Spectral sequence with local cohomology groups}

In this subsection, we construct an alternative spectral sequence 
converging to the Dolbeault cohomology groups $H^*(X,\OO(E))$.

For any locally closed subset $A\subset X$, let $\Gam_A$ be the functor
which associates every sheaf $\FF$ an Abelian group
$\Gam_A(\FF)=\set{s\in\Gam(\FF)}{\supp s\subset A}$.
Recall that the {\em local cohomology groups} $H^q_A$ ($q\ge0$) are 
the derived functors of $\Gam_A$, i.e., $H^q_A(X,\FF)=H^q(\Gam_A(\CF))$.
The {\em sheaves of local cohomology} $\HH^q_A(\FF)$ with supports in $A$
are the sheaves associated to the presheaves $U\mapsto H^q_{U\cap A}(U,\FF)$,
where $U$ is any open subset of $X$.
(We refer the reader to [\ref{BS}, chap.\ II] and [\ref{K}, \S\,7-10] 
for details.)
For any closed subset $A'$ of $A$, let 
$\Gam_{A/A'}(\FF)=\Gam_A(\FF)/\Gam_{A'}(\FF)$.
If $\FF$ is flabby, then $\Gam_{A/A'}(\FF)=\Gam_{A-A'}(\FF)$
[\ref{K}, Lemma~7.3].
The derived functors of $\Gam_{A/A'}$ are denoted by $H^q_{A/A'}$.
We have $H^q_{A/A'}(\FF)=H^q_{A-A'}(\FF)$ for any sheaf $\FF$.
Let $\HH^q_{A/A'}(\FF)$ be the sheaves associated to the presheaves 
$U\mapsto H^q_{U\cap A/U\cap A'}(U,\FF)$, where $U$ is any open subset of $X$.

Again, we denote $\OO(E)$ by $\FF$ from now on.
If $A$ is $T$-invariant, then $H^q_A(X,\FF)$ ($q\ge0$) are representations
of $T$.

\begin{lemma}\label{FILT'}
Under Assumption~\ref{WA}, there is a $T$-equivariant spectral sequence with
	\be\label{E1}
E\pq_1=H^{p+q}_{Z_p-Z_{p+1}}(X,\FF)=
\bigoplus_{\al\in F_p}H^{p+q}_{X^C_\al}(X,\FF)
	\ee
that converges $T$-equivariantly to $H^*(X,\FF)$.
\end{lemma}

\proof{From (\ref{closed}), we have a filtration of the cochain complex
	\be
\Gam(\CF)=\Gam_{Z_0}(\CF)\supset\Gam_{Z_1}(\CF)\supset\cdots\supset
\Gam_{Z_m}(\CF)\supset\Gam_{Z_{m+1}}(\CF)=0.
	\ee
This induces a spectral sequence that converging to $H^*(X,\FF)$ with
	\be
E\pq_0=\Gam_{Z_p}(\CC^{p+q}(\FF))/\Gam_{Z_{p+1}}(\CC^{p+q}(\FF))=
\Gam_{Z_p-Z_{p+1}}(\CC^{p+q}(\FF)).
	\ee
Therefore
	\be
E\pq_1=H^{p+q}_{Z_p-Z_{p+1}}(\CF)=H^{p+q}_{Z_p-Z_{p+1}}(X,\FF).
	\ee
(See for example [\ref{Z}, Theorem~1.1]; the proof is included here for
completeness.)
Since $Z_p-Z_{p+1}=\bigcup_{\al\in F_p}X^C_\al$ and
$\overline{X^C_\al}\cap X^C_\beta=\emptyset$ for $\al\ne\beta\in F_p$, 
we have $\Gam_{Z_p-Z_{p+1}}(\CF)=\bigoplus_{\al\in F_p}\Gam_{X^C_\al}(\CF)$.
Hence
	\be
H^*_{Z_p-Z_{p+1}}(X,\FF)=\bigoplus_{\al\in F_p}H^*_{X^C_\al}(X,\FF).
	\ee}

\vspace{1ex}
Similar to the study of cohomology with compact support,
we consider a general fibration $\pi\colon Y\to B$.
Suppose for the time being that $X$ is any topological space containing $Y$ as 
a locally closed subset and that $\FF$ is any sheaf on $X$.
We want to compute the local cohomology groups $H^q_Y(X,\FF)$ ($q\ge0$).
Let $\A$, $\LL^*$ be the sheaves on $B$ defined by the presheaves
$\A(U)=\Gam_{\pi^{-1}(U)}(X,\FF)$, $\LL^*(U)=\Gam_{\pi^{-1}(U)}(X,\CF)$,
respectively, where $U$ is any open subset of $B$.
Then $0\to\A\to\LL^*$ is a differential sheaf in the sence of
[\ref{G}, \S\,II.4.1].
Let $\HH^q_Y(X,\FF)$ ($q\ge0$) be the sheaves on $B$ 
defined by the presheaves $\HH^q_Y(X,\FF)(U)=H^q(\LL^*(U))$, for any open
subset $U\subset B$.

\begin{lemma}\label{PROJ'}
1. At $b\in B$, the stalk of $\HH^q_Y(X,\FF)$ for any $q\ge0$ is
	\be
\HH^q_Y(X,\FF)_b\cong H^q_{Y_b}(X,\FF).
	\ee
2. There is a spectral sequence with
	\be
E\pq_2=H^p(B,\HH^q_Y(X,\FF))
	\ee
that converges to $H^*_Y(X,\FF)$.
\end{lemma}

\proof{1. This is the analog of Lemma~\ref{PROJ}.1. We have
$\HH^q_Y(X,\FF)_b=\indlim_{U\ni b}H^q_{\pi^{-1}(U)}(X,\FF)=H^q_{Y_b}(X,\FF)$;
the second equality follows from $\bigcap_{U\ni b}\pi^{-1}(U)=Y_b$.\\
2. It is clear that $\LL^*$ are flabby sheaves and that 
$\Gam(B,\LL^*)=\Gam_Y(X,\CF)$.
The rest of the proof is identical to that of Lemma~\ref{PROJ}.2.}

\begin{lemma}\label{MODEL'}
Under the conditions of Lemma~\ref{MODEL}, we have
	\be
H^q_A(R,\OO(E))=
\two{S(A^*)\otimes S(A^\perp)\otimes\medwedge^{n-\nu}(A^\perp)\otimes E_0}
{q=n-\nu}{0}{q\ne n-\nu.}
	\ee
\end{lemma}

\proof{As in the proof of Lemma~\ref{MODEL}, the general result follows from
	\be
H^q_\zero(R,\OO)=\two{S(R)\otimes\medwedge^n(R)}{q=n}{0}{q\ne n}
	\ee
and
	\be
H^q(R,\OO)=\two{S(R^*)}{q=0}{0}{q\ne0.}
	\ee
See also [\ref{K}, Proposition~11.9(e)].}

\vspace{1ex}
We now return to the situation of $\FF=\OO(E)$.

\begin{lemma}\label{EACH'}
	\be\label{each'}
H^q_{X^C_\al}(X,\FF)=H^{q+\nu^C_\al+n_\al-n}(X^T_\al,\OO(S((N^C_\al)^*)
\otimes S(N^{-C}_\al)\otimes\medwedge^{n-n_\al-\nu^C_\al}(N^{-C}_\al)
\otimes E|_{X^T_\al}))
	\ee
as representations of $T$.
\end{lemma}

\proof{Consider the fibration $\pi^C\colon X^C_\al\to X^T_\al$.
For any $x\in X^T_\al$, we want to find the stalk $\HH^q_{X^C_\al}(X,\FF)_x$,
which by excision [\ref{BS}, \S\,II.1, Lemma~1.1] depends only on 
an open neighborhood of $(\pi^C)^{-1}(x)\subset X^C_\al$ in $X$.
By Lemma~\ref{LOCAL}, we can replace $X^C_\al\subset X$ by 
$N^C_\al|_{U_x}\subset N_\al|_{U_x}$ and $E$ by a trivial vector bundle
with fiber $E_x$.
Moreover there is a $T$-equivariant isomorphism
$(N_\al,N^C_\al)|_{U_x}\cong U_x\times(N_x,N^C_x)$.
By Lemma~\ref{PROJ'}.1 and Lemma~\ref{MODEL'},
	\bea
\HH^q_{N^C}(X,\FF)_x\eq H^q_{N^C_x}(W^C_x,\OO(W^C_x,E_x))		\nno
\eq\two{\OO(S((N^C_\al)^*)\otimes S(N^{-C}_\al)\otimes
\medwedge^{n-n_\al-\nu^C_\al}(N^{-C}_\al)\otimes E|_{X^T_\al})}
{q=n-n_\al-\nu^C_\al}{0}{q\ne n-n_\al-\nu^C_\al.}
	\eea
So the spectral sequence of Lemma~\ref{PROJ'}.2 degenerates at $E_2$
and the result follows.}

\begin{thm}\label{MAIN'}
Under the conditions of Theorem~\ref{MAIN},\\
1. there is a $T$-equivariant spectral sequence converging $T$-equivariantly
to $H^*(X,\OO(E))$ with
	\be\label{e1'}
E\pq_1=
\bigoplus_{\al\in F_p}H^{p+q+\nu^C_\al+n_\al-n}(X^T_\al,\OO(S((N^C_\al)^*)
\otimes S(N^{-C}_\al)\otimes\medwedge^{n-n_\al-\nu^C_\al}(N^{-C}_\al)
\otimes E|_{X^T_\al}));
	\ee
2. there is a character valued polynomial $Q^C(t)\ge0$ such that
	\bea\label{morse'}
\vc\sum_{\al\in F}t^{n-n_\al-\nu^C_\al}\sum_{q=0}^{n_\al}t^q
\ch H^q(X^T_\al,\OO(S((N^C_\al)^*)\otimes S(N^{-C}_\al)\otimes
\medwedge^{n-n_\al-\nu^C_\al}(N^{-C}_\al)\otimes E|_{X^T_\al}))		\nno
\eq\sum_{q=0}^nt^q\ch H^q(X,\OO(E))+(1+t)Q^C(t);
	\eea
3. 
	\be\label{ind'}
\sum_{q=0}^n(-1)^q\ch H^q(X,\OO(E))=
\sum_{\al\in F}(-1)^{n-n_\al-\nu^C_\al}\int_{X^T_\al}
{\mathrm{ch}}_T\left(\frac{E|_{X^T_\al}\otimes\det(N^{-C}_\al)}
{\det(1-(N^C_\al)^*)\otimes\det(1-N^{-C}_\al)}\right){\mathrm{td}}(X^T_\al).
	\ee
\end{thm}

\proof{Part~1 follows from Lemma~\ref{FILT'} and Lemma~\ref{EACH'}.
Parts~2 and 3 are proved in the same way as in Theorem~\ref{MAIN}.}

\begin{cor}\label{DISCR'}
Under the conditions of Corollary~\ref{DISCR},\\
1. there is a $T$-equivariant spectral sequence converging $T$-equivariantly
to $H^*(X,\OO(E))$ with
	\be\label{e1-discr'}
E\pq_1=\bigoplus_{x\in F,\,\nu^C_x=n-p-q}S((N^C_x)^*)\otimes S(N^{-C}_x)
\otimes\medwedge^{n-\nu^C_x}(N^{-C}_x)\otimes E_x.
	\ee
2. there is a character valued polynomial $Q^C(t)\ge0$ such that
	\be\label{morse-discr'}
\sum_{x\in F}t^{n-\nu^C_x}\ch(E_x)
\prod_{\lam_{x,k}\in C^*}\inv{1-\e{-\lam_{x,k}}}
\prod_{\lam_{x,k}\in-C^*}\frac{\e{\lam_{x,k}}}{1-\e{\lam_{x,k}}}
=\sum_{q=0}^nt^q\ch H^q(X,\OO(E))+(1+t)Q^C(t);
	\ee
3. 
	\be\label{ind-discr'}
\sum_{q=0}^n(-1)^q\ch H^q(X,\OO(E))=\sum_{x\in F}(-1)^{n-\nu^C_\al}
\ch(E_x)\prod_{\lam_{x,k}\in C^*}\inv{1-\e{-\lam_{x,k}}}
\prod_{\lam_{x,k}\in-C^*}\frac{\e{\lam_{x,k}}}{1-\e{\lam_{x,k}}}.
	\ee
\end{cor}

\begin{rmk}\label{REMARK'}{\em
1. The same observations in Remark \ref{REMARK}.1 apply to Theorem~\ref{MAIN'}
and Corollary~\ref{DISCR'}.
In particular, all the weights of $T$ in $H^*(X,\OO(E))$ are also of finite
multiplicities.
When $X$ is non-compact, the Dolbeault cohomology groups are different from 
those with compact support.
Therefore the results of Theorem~\ref{MAIN} and Theorem~\ref{MAIN'} are 
not the same.
Again, it would be interesting to have an independent analytic proof of 
parts~2 of Theorem~\ref{MAIN'} and Corollary~\ref{DISCR'} when $X$ is 
a non-compact K\"ahler manifold satisfying Assumption~\ref{PWA}.
When $X$ is compact, Theorem~\ref{MAIN} is identical to Theorem~\ref{MAIN'}
with an opposite action chamber.
It is possible that the two theorems are dual to each other in some sense;
this is also reflected by the local models in Lemma~\ref{MODEL} and 
Lemma~\ref{MODEL'}.\\
2. Remark \ref{REMARK}.2 applies here as well.
In particular, the complex $\{E^{*0}_1,\dr^{*0}_1\}$, i.e.,
	\be\label{gC}
0\to\Gam(X,\FF)\to H^0_{Z_0-Z_1}(X,\FF)\to H^1_{Z_1-Z_2}(X,\FF)\to\cdots
\to H^m_{Z_m}(X,\FF)\to0,
	\ee
is called the {\em global Grothendieck-Cousin complex} [\ref{H}, \ref{K}].
If condition (\ref{dege2}) is satisfied, then 
the complex (\ref{gC}) computes the cohomology groups $H^*(X,\OO(E))$.
Again a sufficient condition for (\ref{dege2}) is that $X^T=F$ is discrete, 
$m=n$ in (\ref{closed}), and $F_p=\set{x\in F}{\nu^{-C}_x=p}$ for 
$0\le p\le n$.
In [\ref{K}, \S\,10], a few other sufficient conditions were found.
If $\HH^q_{Z_p/Z_{p+1}}(\FF)=0$ for all $q\ne p$, then the complex of sheaves
	\be\label{lC}
0\to\FF\to\HH^0_{Z_0/Z_1}(\FF)\to\HH^1_{Z_1/Z_2}(\FF)\to\cdots\to
\HH^m_{Z_m}(\FF)\to0,
	\ee
called the {\em local Grothendieck-Cousin complex}, is a resolution of $\FF$
(see for example [\ref{K}, Theorem~8.7] or [\ref{Boz}, Lemma~1.2]).
In this case, the sheaf $\FF$ is called {\em locally Cohen-Macaulay}
with respect to the filtration (\ref{closed}).
The global Grothendieck-Cousin complex (\ref{gC}),
which computes the cohomology groups $H^*(X,\FF)$,
is obtained from (\ref{lC}) by applying the functor $\Gam(X,\cdot\,)$.}
\end{rmk}

\sect{Examples and Applications}

\subsect{Calculations in \v Cech coholomogy theory}

We interpret some of the procedures in the last section in the language of
\v Cech cohomology theory, which is especially suitable for calculations.

For any sheaf $\FF$ on $X$ and any open cover $\UU=\set{U_i}{i\in I}$ of $X$,
let $H^*(\UU,\FF)$ be the cohomology groups of the \v Cech cochain complex
	\be
C^q(\UU,\FF)=\bigoplus_\iii\FF(\UUU),
	\ee
where $\iii\in I$ ($q\ge0$) are not equal, with the standard coboundary maps.
The \v Cech cohomology groups $\check H^*(X,\FF)$ are the inductive limits 
of $H^*(\UU,\FF)$ with respect to the refinement of open coverings.
We have the well-known isomorphism $\check H^*(X,\FF)\cong H^*(X,\FF)$.
For any open subset $V\subset X$, let $\FFt_V$ be the presheaf defined by
	\be
\FFt_V(U)=\two{\FF(U)}{U\subset V}{0}{\mbox{otherwise}.}
	\ee
Then $\FF_V$ is the sheaf associated to the presheaf $\FFt_V$ and
$\check H^*(X,\FFt_V)=\check H^*(X,\FF_V)$ [\ref{G}, \S II.5.11].
If $V'$ is an open subset of $V$, then the presheaf $\FFt_V/\FFt_{V'}$,
denoted by $\FFt_{V/V'}$, generates the sheaf $\FF_V/\FF_{V'}=\FF_{V-V'}$.
Moreover, $\check H^*(X,\FFt_{V/V'})=\check H^*(X,\FF_{V-V'})$.

We now assume that $X$ is a compact complex manifold with a meromorphic
$T$-action and that there is an action chamber $C$ such that the set $F$ of 
connected components of $X^T$ is partially ordered with respect to 
the relation $\prec$.
Then the \bb decomposition is filterable, with filtrations of $X$ by
closed subsets (\ref{closed}) and by open subsets (\ref{open}).
Let $\FF=\OO(E)$, where $E\to X$ is a holomorphic vector bundle on which the
$T$-action lifts holomorphically.
Then there is a filtration of $\FF$ by presheaves
	\be
\FF=\FFt_{V_0}\supset\FFt_{V_1}\supset\cdots\supset\FFt_{V_m}
\supset\FFt_{V_{m+1}}=0.
	\ee
We choose $\UU$ to be a $T$-invariant, i.e., for any $U_i\in\UU$, $g\in T$,
we have $gU_i\in\UU$.
Then we have a $T$-equivariant filtration of the \v Cech complex 
	\be\label{fcech}
C^*(\UU,\FF)=F^0C^*\supset F^1C^*\supset\cdots\supset F^mC^*\supset
F^{m+1}C^*=0,
	\ee
where
	\be\label{fcechpq}
F^pC^q=\bigoplus_\iii\FFt_{V_p}(\UUU)=\bigoplus_{\UUU\subset V_p}\FF(\UUU).
	\ee
So there is a $T$-equivariant spectral sequence converging to $H^*(\UU,\FF)$
with $E\pq_0=\bigoplus_\iii\FFt_{V_p/V_{p+1}}(\UUU)$ and 
$E\pq_1=H^{p+q}(\UU,\FFt_{V_p/V_{p+1}})$.
Taking the inductive limit of $\UU$, we conclude that there is 
a $T$-equivariant spectral sequence with $E\pq_1=\check H(X,\FF_{V_p-V_{p+1}})$
that converges to $\check H^*(X,\FF)$.
This is Lemma~\ref{FILT} when $X$ is compact.
However the method outlined above is more convenient for calculations, which 
we now illustrate.

\begin{ex}\label{CP2}{\em
Let $X=\CP^2=\{[z_0,z_1,z_2]\}$ be equipped with a holomorphic action of
$T^2=\CM\times\CM$ given by
$(u_1,u_2)\colon[z_0,z_1,z_2]\mapsto[z_0,u_1^{-1}z_1,u_2^{-1}z_2]$.
Let $\{e_1,e_2\}$ be the standard basis of $\gt_\re\cong\RE^2$, 
and $\{e_1^*,e_2^*\}$, the dual basis of $\gt_\re^*\cong\RE^2$.
The three isolated fixed points of $T^2$ in $X$ are $p_0=[1,0,0]$, 
$p_1=[0,1,0]$, $p_2=[0,0,1]$, whose isotropy weights are
$-e_1^*$, $-e_2^*$; $e_1^*$, $e_1^*-e_2^*$; $e_2^*$, $e_2^*-e_1^*$.
The $T^2$-action on $X=\CP^2$ is meromorphic.
Choose the action chamber $C$ spanned by $e_1+e_2$ and $e_2$.
Then the relation $\prec$ on $F=\{p_0,p_1,p_2\}$ is given by
$p_2\prec p_1\prec p_0$.
The cells in the \bb decomposition (with respect to $C$) are
$X^C_0=\{[1,0,0]\}=\{p_0\}$, $X^C_1=\{[z_0,1,0]\}\cong\CO$, 
$X^C_2=\{[z_0,z_1,1]\}\cong\CO^2$.
$X$ has a filtration by closed subsets 
$X=Z_0\supset Z_1\supset Z_2\supset Z_3=\emptyset$
with $Z_1=\{[z_0,z_1,0]\}\cong\CP^1$, $Z_2=\{[1,0,0]\}=\{p_0\}$
and a filtration by open subsets
$X=Z_0\supset V_1\supset V_2\supset V_3=\emptyset$ with
$V_1=\set{[z_0,z_1,z_2]}{z_1\ne0\mbox{ or }z_2\ne0}$, $V_2=\{[z_0,z_1,1]\}$.
We have $Z_{2-q}-Z_{3-q}=V_q-V_{q+1}=X_q^C$ ($0\le q\le2$).
Let $L=((\CO^3-\zero)\times\CO)/\CM$, where the $\CM$-action is
$u\colon(z_0,z_1,z_2,w)\mapsto(uz_0,uz_1,uz_2,u^cw)$ for some $c\in\ZZ$.
$L$ is a holomorphic line bundle over $X=\CP^2$ with a lifted holomorphic 
$T^2$-action
$(u_1,u_2)\colon[z_0,z_1,z_2,w]\mapsto[z_0,u_1^{-1}z_1,u_2^{-1}z_2,w]$.
The weights on the fibers $L_0$, $L_1$, $L_2$ over $p_0$, $p_1$, $p_2$
are $0$, $ce_1^*$, $ce_2^*$, respectively.
The first Chern class of $L$ is $c_1(L)=c$.
Let $\FF=\OO(L)$.
We want to calculate $H^*(X,\FF)$ in \v Cech theory.

Choose an open covering $\UU=\{U_0,U_1,U_2\}$, where
$U_i=\set{[z_0,z_1,z_2]}{z_i\ne0}\cong\CO^2$ ($i=0,1,2$).
Then $U_i\cap U_j\cong\CO\times\CM$ for any $0\le i<j\le2$
and $U_0\cap U_1\cap U_2\cong\CM\times\CM$.
The restrictions of $L$ to $U_i$ are trivial, and
$\FF(U_i)=\OO(U_1)\otimes L_i$ ($i=0,1,2$).
Since $H^1(\CO,\OO)=H^1(\CM,\OO)=0$, the open cover $\UU$ already satisfies
$H^*(\UU,\FF)\cong\check H^*(X,\FF)$.
Therefore $\check H^*(X,\FF)$ can be computed by the spectral sequence 
associated to the filtration (\ref{fcech}).
According to (\ref{fcechpq}), the spaces $F^pC^q$ are given by
	\begin{center}
	\begin{tabular}{c|cccc}
$q=2$ &	$\FF(U_0\cap U_1\cap U_2)$	&  $\FF(U_0\cap U_1\cap U_2)$	&
	$\FF(U_0\cap U_1\cap U_2)$	&	$0$			\\
$q=1$ &	$\bigoplus_{i<j}\FF(U_i\cap U_j)$&  $\bigoplus_{i<j}\FF(U_i\cap U_j)$ &
	$\FF(U_0\cap U_2)\oplus\FF(U_1\cap U_2)$  &	$0$		\\
$q=0$ &	$\bigoplus_i\FF(U_i)$		&  $\FF(U_1)\oplus\FF(U_2)$	&
	$\FF(U_2)$			&	$0$			\\
\cline{2-5}
\multicolumn{1}{c}{} & 	$p=0$		&	$p=1$			&
	$p=2$				&	$p=3$
	\end{tabular}
	\end{center}
Consequently, the spaces $E\pq_0=F^pC^{p+q}/F^{p+1}C^{p+q}$ are given by
	\begin{center}
	\begin{tabular}{r|ccccc}
\multicolumn{1}{r}{} & &   $p=0$ & &	$p=1$	&	$p=2$		\\
\cline{2-6}
$q=0$  & & $\FF(U_0)$ & & $\FF(U_0\cap U_1)$	&  $\FF(U_0\cap U_1\cap U_2)$\\
\vspace{-5ex}								\\
       & &	      &	& $\uparrow$	&	$\uparrow$		\\
\vspace{-5ex}								\\
$q=-1$ & & & & $\FF(U_1)$ &  $\FF(U_0\cap U_2)\oplus\FF(U_1\cap U_2)$	\\
\vspace{-5ex}								\\
       & &  &	      &	&	$\uparrow$				\\
\vspace{-5ex}								\\
$q=-2$ & &  &	      &	&  $\FF(U_2)$
	\end{tabular}
	\end{center}
The arrows here (and below) denote the coboundary operators.
Let $\Gam\pqr$ be a region in $\gt^*_\re\cong\RE^2=\{(x,y)\}$ such that
$\supp E\pqr=\Gam\pqr\cap\ell^*$.
(Recall that $\ell^*\cong\ZZ^2$ is the dual lattice in $\gt_\re^*\cong\RE^2$.)
Then regions $\Gam\pq_0$ are given by
	\begin{center}
	\begin{tabular}{r|cccc}
\multicolumn{1}{c}{}  & &   $p=0$   &	$p=1$		&	$p=2$	\\
\cline{2-5}
$q=0$  & & $\{x,y\ge0\}$ 	& 	$\{y\ge0\}$	&  $\{(x,y)\}$	\\
\vspace{-5ex}								\\
       & &		& $\uparrow$	&	$\uparrow$		\\
\vspace{-5ex}								\\
$q=-1$ & &	& $\{y\ge0,\;x+y\le c\}$ &  $\{x\ge0\}\uplus\{x+y\le c\}$\\
\vspace{-5ex}								\\
       & &  &		&	$\uparrow$				\\
\vspace{-5ex}								\\
$q=-2$ & &		  &			&  $\{x\ge0,\;x+y\le c\}$
	\end{tabular}
	\end{center}
Here the multiplicity of any weight $\xi\in\Gam\pq_0\cap\ell^*$ in $E\pq_0$
is $1$ except for $E^{2,-1}_0$ and $\xi\in\{x\ge0\}\cap\{x+y\le c\}$, in
which case the multiplicity is $2$.
The regions $\Gam\pq_1$ are given by
	\begin{center}
	\begin{tabular}{r|cccccc}
\multicolumn{1}{c}{}  & &  $p=0$   &	& $p=1$	    &       & $p=2$	\\
\cline{2-7}
$q=0$	& & $\{x,y\ge0\}$ & $\to$ & $\{y\ge0,\;x+y>c\}$ & $\to$ & 
$\{x<0,\,x+y<c\}$
	\end{tabular}
	\end{center}
Finally, the regions $\Gam\pq_2$ are given by
	\begin{center}
	\begin{tabular}{r|ccccccccc}
\multicolumn{1}{c}{} &   $p=0$   &	$p=1$	   & &	& &	$p=2$ & & & \\
\cline{2-10}
$q=0$	& $\{x,y\ge0,\;x+y\le c\}$ & &	& &	&	& &	&
	\end{tabular}
	\end{center}
if $c\ge0$ and by
	\begin{center}
	\begin{tabular}{r|ccccccccc}
\multicolumn{1}{c}{} & & & &  $p=0$	& & &	 &  $p=1$	&  $p=2$  \\
\cline{2-10}
$q=0$	& & & &		& & &	 &		& $\{x,y<0,\;x+y>c\}$
	\end{tabular}
	\end{center}
if $c\le-2$.
All the weights $\xi\in\Gam\pqr\cap\ell^*$ ($r=1,2$) are of multiplicity $1$ 
in $E\pqr$.
So the spectral sequence degenerates at $E_2$.
We recover the well-known result that the only non-trivial cohomology groups 
are $H^0(\CP^2,\OO(L))$ if $c\ge0$ and $H^2(\CP^2,\OO(L))$ if $c\le-2$
for a holomorphic line bundle $L\to\CP^2$ of $c_1(L)=c$.}
\end{ex}

\begin{rmk}{\em 
The method of Example~\ref{CP2} applies to any toric variety satisfying
the partial order condition. 
More interestingly, the (holomorphic) instanton complex can be used to study
the cohomology groups of vector bundles over spherical varieties, about
which not all is known. 
(See [\ref{Br}] for an extension of the Borel-Weil theorem.)
One notable exception is the flag manifold, which will be discussed
in the next subsection.}
\end{rmk}

\subsect{Flag manifolds and generalized Bernstein-Gelfand-Gelfand resolutions}

We show that the spectral sequence for the cohomology of a flag manifold
leads to geometric realizations of the Bernstein-Gelfand-Gelfand [\ref{BGG}]
and related resolutions.

Let $G$ be a complex semi-simple Lie group and $T$, a maximal torus of $G$.
Let $\g$, $\gt$ be the Lie algebras of $G$, $T$, respectively.
Let $\g=\gt\oplus\bigoplus_{\al\in\Del}\g_\al$ be the root space decomposition,
where $\Del\subset\gt^*-\zero$ is the root system of the pair $(\g,\gt)$
and $\g_\al=\CO e_\al$ ($\al\in\Del$).
Let $\Del_+$ be a set of positive roots and let $\Del_-=-\Del_+$.
Let $\gn_\pm=\bigoplus_{\al\in\Del_\pm}\g_\al$.
Let $B$ be the Borel subgroup corresponding to the Borel subalgebra
$\gb=\gt\oplus\gn_+$.
Let $W$ be the Weyl group of the pair $(\g,\gt)$.
Denote by $w_0$ the element in $W$ of maximal length $l(w_0)=|\Del_+|$.

Recall that the Verma module of highest weight $\lam$ is the $U(\g)$-module
$M_\lam=U(\g)\otimes_{U(\sgb)}\CO v_\lam$, where $\CO v_\lam$ is the
$1$-dimensional $U(\gb)$-module defined by $\lam\in\gt_\re^*$.
$M_\lam$ is free over $U(\gn_-)$.
As a $U(\gt)$-module, $M_\lam$ is determined by
$\ch M_\lam=\frac{\e{\lam}}{\prod_{\al\in\Del_+}(1-\e{-\al})}$.
When $\lam$ is a dominant weight, let $R_\lam$ be the (finite dimensional)
irreducible module of highest weight $\lam$.
We have a resolution of $R_\lam$ by Verma modules [\ref{BGG}]
	\be\label{bgg}
0\to M_{w_0\lam-2\rho}\to\bigoplus_{l(w)=|\Del_+|-1}M_{w(\lam+\rho)-\rho}
\to\cdots\to\bigoplus_{l(w)=1}M_{w(\lam+\rho)-\rho}\to M_\lam\to R_\lam\to0,
	\ee
where $\rho=\hf\sum_{\al\in\Del_+}\al$.
This is called the Bernstein-Gelfand-Gelfand resolution of $R_\lam$.

For any $w\in W$, put $\gn^w_\pm=w\gn_\pm w^{-1}$.
The {\em twisted Verma module} $M^w_\lam$ is a $U(\g)$-module of highest 
weight $\lam$ that is free over $U(\gn^w_+\cap\gn_-)$ and co-free over
$U(\gn^w_+\cap\gn_+)$ [\ref{FF}].
In particular, $M^1_\lam\cong M^*_\lam$ and $M^{w_0}_\lam\cong M_\lam$;
the $U(\g)$-module structure of the dual $M^*_\lam$ is given by
[\ref{BK}, \S 2.3]
	\be
\bra x\xi,v\ket=-\bra\xi,\tau(x)v\ket\quad\mbox{for }x\in\g,\;\xi\in M_\lam^*,
\;v\in M_\lam,
	\ee
where $\tau$ is an automorphism of $\g$ such that $\tau(h)=-h$ ($h\in\gt$)
and $\tau(e_\al)=e_{-\al}$ ($\al\in\Del$).
If $M_\lam$ is irreducible, then $M^w_\lam\cong M_\lam$ for any $w\in W$.
As $U(\gt)$-modules, we always have $\ch M^w_\lam=\ch M_\lam$.

For example, take $\g=\gsl(2,\CO)=\spanc\{h,e,f\}$ with commutation
relations
	\be
[h,e]=2e,\quad [h,f]=-2f,\quad [e,f]=h.
	\ee
The Weyl group is $W=\{\pm1\}$.
The twisted Verma modules of highest weight $\lam\in\RE$ are
$M^1_\lam=M^*_\lam$ and $M^{-1}_\lam=M_\lam$.
Here the Verma module is $M_\lam=\spanc\set{v_\lam^k=f^kv_\lam}{k\in\NN}$ with
	\be
\three{hv_\lam^k=(\lam-2k)v_\lam^k}{ev_\lam^k=k(\lam-k+1)v_\lam^{k-1}}
{fv_\lam^k=v_\lam^{k+1}}\quad(k\in\NN).
	\ee
The dual is $M_\lam^*=\spanc\set{\xi_\lam^k}{k\in\NN}$, where
$\bra\xi_\lam^k,v_\lam^l\ket=\del^{kl}$ ($k,l\in\NN$).
The $U(\g)$-module structure is
	\be
\three{h\xi_\lam^k=(\lam-k)\xi_\lam^k}{e\xi_\lam^k=-\xi_\lam^{k-1}}
{f\xi_\lam^k=-(k+1)(\lam-k)\xi_\lam^{k+1}}\quad(k\in\NN).
	\ee
When $\lam$ is not a dominant weight, i.e., when $\lam\not\in\NN$,
$M_\lam\cong M_\lam^*$ as $U(\g)$-modules and the isomorphism is given by
$v_\lam^k\mapsto k!(k-1-\lam)\cdots(-\lam)\xi_\lam^k$ ($k\in\NN$).
When $\lam\in\NN$, $M_\lam$ and $M_\lam^*$ are not isomorphic $U(\g)$-modules.
The irreducible module $R_\lam$ is a quotient of $M_\lam$ and a submodule of
$M_\lam^*$.

We consider the non-degenerate flag manifold $X=G/B^-$, where $B^-$ is the 
Borel subgroup opposite to $B$.
The maximal torus $T$ acts meromorphically on $X$.
The fixed-point set is $X^T=\set{wB^-}{w\in W}$.
The isotropy weights at $wB^-$ are $w\al$ ($\al\in\Del_+$).
The action chambers  in $\gt$ are the Weyl chambers.
Choose the positive Weyl chamber, denoted by ``+''.
Then the polarizing index of $wB^-$ is $\nu^-_w=|\Del_-\cap w\Del_+|=l(w)$
for any $w\in W$.
The \bb decomposition is precisely the Bruhat decomposition 
$X=\bigcup_{w\in W}X^+_w$, where $X^+_w=BwB^-/B^-$ ($w\in W$) are the Bruhat
cells [\ref{Ak}].
These cells are also the $B$-orbits in $X$.
Moreover, the relation $\prec$ on $F\cong W$ is the Chevalley-Bruhat
order [\ref{C}], which is a partial ordering.
Consequently, the \bb decomposition is filterable, and we have the filtration
(\ref{closed}), where $m=|\Del_+|=\dim_\co X$.
The closed sets $Z_p=\bigcup_{l(w)\ge p}X^+_w$ ($0\le p\le|\Del_+|$) are 
the Schubert varieties.
Since $Z_p-Z_{p+1}=\bigcup_{l(w)=p}X^+_w$ ($0\le p\le|\Del_+|$) and 
$\nu^-_w=l(w)$, 
the cohomology groups $H^*(X,\FF)$ with coefficients in any sheaf $\FF$ can be
computed by the (global) Grothendieck-Cousin complex (\ref{gC}), which becomes
	\be\label{kempf}
0\to H^0_{X^+_1}(X,\FF)\to\bigoplus_{l(w)=1}H^1_{X^+_w}(X,\FF)\to\cdots\to
\bigoplus_{l(w)=|\Del_+|-1}H^{|\Del_+|-1}_{X^+_w}(X,\FF)\to 
H^{|\Del_+|}_{X^+_{w_0}}(X,\FF)\to0.
	\ee

Given any integral weight $\lam\in\ell^*$, we have a holomorphic line bundle
$L_\lam=G\times_{B^-}\CO v_\lam$ over $X$, where $\CO v_\lam$ is the
$1$-dimensional holomorphic representation of $B^-$ defined by $\lam$.
The weight of $T$ on the fiber $(L_\lam)_{wB^-}$ ($w\in W$) is $w\lam$.
Set $\FF_\lam=\OO(L_\lam)$.
Then from subsection~3.2, we have for any $w\in W$,
	\be
\ch H^{l(w)}_{X^+_w}(X,\FF_\lam)=\e{w\lam}
\prod_{\al\in\Del_+\cap w^{-1}\Del_+}\inv{1-\e{-w\al}}
\prod_{\al\in\Del_+\cap w^{-1}\Del_-}\frac{\e{w\al}}{1-\e{w\al}}
=\frac{\e{w(\lam+\rho)-\rho}}{\prod_{\al\in\Del_+}(1-\e{-\al})}.
	\ee
So as representations of $T$, $H^{l(w)}_{X^+_w}(X,\FF_\lam)$ is the same as the
Verma module $M_{w(\lam+\rho)-\rho}$.
In fact the above local cohomology groups are $U(\g)$-modules.
This is because the canonical resolution $\CC^*(\FF_\lam)$ of $\FF_\lam$,
on which the Lie algebra $\g$ acts, is $U(\g)$-equivariant.
(Notice however that a representation of $\g$ on an infinite dimensional
space may not exponentiate to that of $G$.)
Therefore the Grothendieck-Cousin complex (\ref{kempf}) is $U(\g)$-equivariant.
Moreover, we have $H^{l(w)}_{X^+_w}(X,\FF_\lam)\cong M^w_{w(\lam+\rho)-\rho}$
as $U(\g)$-modules [\ref{FF}, \S 2.2].
So (\ref{kempf}) becomes
	\be\label{gbgg}
0\to M^*_\lam\to\bigoplus_{l(w)=1}M^w_{w(\lam+\rho)-\rho}\to\cdots\to
\bigoplus_{l(w)=|\Del_+|-1}M^w_{w(\lam+\rho)-\rho}\to M_{w_0\lam-2\rho}\to0.
	\ee
If $\lam+\rho$ is regular, the cohomology groups $H^*(X,\FF_\lam)$,
hence those of the complex (\ref{gbgg}), are
	\be
H^q(X,\FF_\lam)=
\two{R_{w_\lam(\lam+\rho)-\rho}}{q=l(w_\lam)}{0}{q\ne l(w_\lam),}
	\ee
where $w_\lam$ is the unique element in $W$ such that $w_\lam(\lam+\rho)-\rho$
is a dominant weight [\ref{B1}].
The complex (\ref{gbgg}) is called the {\em generalized 
Bernstein-Gelfand-Gelfand resolution} of
$R_{w_\lam(\lam+\rho)-\rho}$ [\ref{FF}, \S 2.3].
If $\lam+\rho$ is singular, then all $H^*(X,\FF_\lam)=0$ [\ref{B1}].
When $\lam$ is a dominant weight, (\ref{gbgg}) is the dual of the
Bernstein-Gelfand-Gelfand resolution (\ref{bgg}) for $R_\lam$
[\ref{K}, \ref{B}, \ref{BK}].
When $w_0\lam-2\rho$ is dominant, (\ref{gbgg}) is the 
Bernstein-Gelfand-Gelfand resolution for $R_{w_0\lam-2\rho}$.

\begin{ex}{\em
Let $G=SL(2,\CO)$.
Then $T=\set{{u\quad\quad\choose\quad u^{-1}}}{u\in\CM}$ and
$B^-=\set{{u\quad\quad\choose \;*\!\quad u^{-1}}}{u\in\CM}$.
Consider $X=G/B^-=\CP^1=(\CO^2-\zero)/\CM$.
The action of $G$ on the homogeneous coordinates is
	\be
{a\quad b\choose c\quad d}\colon[z_0,z_1]\mapsto[cz_1+dz_0,az_1+bz_0].
	\ee
On the open dense subset $\{[1,z]\}\cong\CO$, the above action is
the fractional linear transformation $z\mapsto\frac{az+b}{cz+d}$.
The generators of $\g$ act on $X$ as $h=2z\ddz$, $e=\ddz$, $f=-z^2\ddz$.
The fixed-point set is $X^T=\{0=[1,0],\,\infty=[0,1]\}$, on which
the Weyl group $W=\{\pm1\}$ acts.
The \bb decompositions $X=X_0^\pm\cup X_\infty^\pm$ were discussed in
subsection~2.1.
For $\lam\in\ZZ$, consider the line bundle is $L_\lam=(\CO^2-\zero)/\CM$,
where the $\CM$-action is $u\colon(z_0,z_1,w)\mapsto(uz_0,uz_1,u^\lam w)$.
The $G$-action lifts to $L_\lam$ according to
	\be
{a\quad b\choose c\quad d}\colon[z_0,z_1,w]\mapsto[cz_1+dz_0,az_1+bz_0,w].
	\ee
Let $\FF_\lam=\OO(L_\lam)$.
The cohomology groups $H^*(\CP^1,\FF_\lam)$ as $U(\g)$-modules can be
computed from the $U(\g)$-equivariant cochain complex
	\be\label{sl2}
0\to H^0_{X^+_0}(\CP^1,\FF_\lam)\to H^1_{X^+_\infty}(\CP^1,\FF_\lam)\to0.
	\ee
Let $U_i=\set{[z_0,z_1]}{z_i\ne0}$ ($i=0,1$) be two open sets in $\CP^1$.
Then the above cochain complex becomes
	\be
0\to\Gam(U_0,\FF_\lam)\to\Gam(U_0\cap U_1,\FF_\lam)/\Gam(U_1,\FF_\lam)\to0.
	\ee
We define two sections $s_i\in\Gam(U_i,\FF_\lam)$ ($i=0,1$) by
$s_i([z_0,z_1])=[z_0,z_1,z_i^\lam]$.
Then $\Gam(U_0,\FF_\lam)=\spanc\set{z^ks_0}{k\in\NN}$ and
$\Gam(U_0\cap U_1,\FF_\lam)/\Gam(U_1,\FF_\lam)=
\spanc\set{z^{k+1}s_1}{k\in\NN}$, where $z=\frac{z_1}{z_0}$ on $U_0$.
The actions of $\g$ on the two spaces are given by
	\be
\three{h(z^ks_0)=(\lam-2k)z^ks_0}{e(z^ks_0)=-kz^{k-1}s_0}
{f(z^ks_0)=-(\lam-k)z^{k+1}s_0}\quad\mbox{and}\quad
\three{h(z^{k+1}s_1)=(-\lam-2k-2)z^{k+1}s_1}
{e(z^{k+1}s_1)=(-\lam-k-1)z^ks_1}{f(z^{k+1}s_1)=(k+1)z^{k+2}s_1.}
	\ee
Therefore as $U(\g)$-modules, $\Gam(U_0,\FF_\lam)\cong M^*_\lam$ and
$\Gam(U_0\cap U_1,\FF_\lam)/\Gam(U_1,\FF_\lam)\cong M_{-\lam-2}$,
where the isomorphisms are given by $z^k\mapsto k!\,\xi^k_\lam$ and
$z^{k+1}s_1\mapsto k!\,v^k_\lam$, respectively.
If $\lam\ge0$, then $M_{-\lam-2}\cong M^*_{-\lam-2}$, and (\ref{sl2}) becomes
$0\to M^*_\lam\to M^*_{-\lam-2}\to0$.
So $H^0(\CP^1,\FF_\lam)=\ker(M^*_\lam\to M^*_{-\lam-2})=R_\lam$
and $H^1(\CP^1,\FF_\lam)=0$.
If $\lam\le-2$, then $M^*_\lam\cong M_\lam$, and (\ref{sl2}) becomes
$0\to M_\lam\to M_{-\lam-2}\to0$.
So $H^0(\CP^1,\FF_\lam))=0$ and 
$H^1(\CP^1,\FF_\lam)=M_{-\lam-2}/M_\lam=R_\lam$.
If $\lam=-1$, then $M^*_\lam\cong M_{-\lam-2}$. 
So all $H^*(\CP^1,\FF_\lam)=0$.}
\end{ex}

\begin{rmk}{\em
Lepowsky [\ref{Lep}] found a Bernstein-Gelfand-Gelfand-type resolution of 
any irreducible $U(\g)$-module by the generalized Verma modules,
which are induced from representations of a parabolic subgroup $P\subset G$.
In [\ref{MR}], a geometric realization of this resolution was constructed
using the local cohomology of the $P$-orbits in $G/B^-$ (rather than
the $B$-orbits in $G/P^-$).
Let $H$ be the Levi subgroup of $P$, and $\gh$, its Lie algebra.
Let $\Del_H$ be the root system of the pair $(\gh,\gt)$, and $W_H$,
the corresponding Weyl group.
Then $X=G/B^-$ decomposes into its $P$-orbits according to
	\be\label{decomP}
X=\bigcup_{w'\in W/W_H}Pw'B^-/B^-,
	\ee
where $W/W_H=\set{W_Hw}{w\in W}$ (see for example [\ref{Wa}, \S 1.2]).
$H$ is the centralizer of a torus subgroup $T'\subset T$, whose Lie algebra
$\gt'\subset\gt$.
Consider the (meromorphic) $T'$-action on $X$.
The fixed-point set $X^{T'}=\bigcup_{w'\in W/W_H}Hw'B^-/B^-$.
Choose the action chamber $C'\subset\gt'$ such that $\bra\al,C'\ket>0$
for all $\al\in\Del_+-\Del_H\cap\Del_+$.
The the \bb decomposition of $X$ with respect to $C'$ is precisely
(\ref{decomP}).
Therefore (\ref{gC}) gives the geometric realizations of Lepowsky's resolution
and similar generalizations.}
\end{rmk}

\subsect{Cohomology and geometric quantization of non-compact manifolds}

In section~3, we obtained equivariant holomorphic Morse inequalities and
equivariant index theorems for non-compact complex manifolds under 
Assumption~\ref{WA}.
In this subsection, we apply them to establish some results on the cohomology
groups and on geometric quantization.

Let $X$ be a (possibly non-compact) complex manifold of dimension $n$
with a holomorphic $T$-action satisfying Assumption~\ref{WA}.
Let $H\pq(X)=H^q(M,\OO(\medwedge^pTX))$,
$H\pq\cpct(X)=H^q\cpct(M,\OO(\medwedge^pTX))$
($p,q=0,1,\dots,n$) be the Dolbeault cohomology groups of $X$ and those
with compact support, respectively.
Let $P(X;s,t)=\sum^n_{p,q=0}s^pt^q\ch H\pq(X)$,
$P\cpct(X;s,t)=\sum^n_{p,q=0}s^pt^q\ch H\pq\cpct(X)$, 
the character-valued Poincar\'e-Hodge polynomials.
If the cohomology groups are finite dimensional, then
$h\pq(X)=\dim_\co H\pq(X)$, $h\pq\cpct(X)=\dim_\co H\pq\cpct(X)$ 
are the Hodge numbers of $X$ and $p(X;s,t)=\sum^n_{p,q=0}s^pt^qh\pq(X)$,
$p\cpct(X;s,t)=\sum^n_{p,q=0}s^pt^qh\pq\cpct(X)$, 
the (usual) Poincar\'e-Hodge polynomials.
Notice that if $X$ is a non-compact \ka manifold, the Hodge numbers or the
Poincar\'e-Hodge polynomials do not necessarily satisfy the usual symmetry
relations.
For example let $X=\CO$.
Then $H^{01}(X)=H^{10}\cpct(X)=0$ whereas $H^{10}(X)$ and $H^{01}\cpct(X)$
are infinite dimensional.

\begin{prop}\label{WZ4.1}
Under Assumption~\ref{WA},\\
1. $\supp H\pq\cpct(X)\subset\overline{C^*}\cap\ell^*$ for all $C$ such that 
the $T$-action is $C$-meromorphic.
Moreover, for any such $C$, there is a polynomial $q^C\cpct(s,t)\ge0$ such that
	\be
\sum_{\al\in F}(st)^{\nu^C_\al}p\cpct(X^T_\al;s,t)=
\sum_{p,q=1}^ns^pt^q\dim_\co H\pq\cpct(X)^T+(1+t)q^C\cpct(s,t);
	\ee
2. $\supp H\pq\cpct(X)\subset-\overline{C^*}\cap\ell^*$ for all $C$ such that
the $T$-action is $C$-meromorphic.
Moreover, for any such $C$, there is a polynomial $q^C(s,t)\ge0$ such that
	\be
\sum_{\al\in F}(st)^{n-n_\al-\nu^C_\al}p(X^T_\al;s,t)=
\sum_{p,q=1}^ns^pt^q\dim_\co H\pq(X)^T+(1+t)q^C(s,t).
	\ee
\end{prop}

\proof{The results follow from the proof of [\ref{WZ}, Theorem 4.1].}

\begin{rmk}\label{WZ4.2}{\em
1. If in addition there is an action chamber $C$ such that the $T$-action 
is both $C$-meromorphic and $(-C)$-meromorphic, then the cohomology groups
$H\pq\cpct(X)$ and $H\pq(X)$ are trivial representations of $T$.
This is true when $X$ is compact [\ref{WZ}, Theorem~4.1.1, Remark 4.2.1]
but not so in general.
For example, let $X=\CO$ with the standard multiplication by $\CM$,
which is plus-meromorphic.
Then $\supp H^{00}(X)=-\NN$ and $\supp H^{01}\cpct(X)=\NN-\zero$.\\
2. As in [\ref{WZ}, Corollary 4.5], we conclude from Proposition~\ref{WZ4.1}
that if $|p-q|>\max_{\al\in F}n_\al$, then $H\pq\cpct(X)^T=H\pq(X)^T=0$.
In particular, if all the fixed points are isolated,
then $H\pq\cpct(X)^T=H\pq(X)^T=0$ when $p\ne q$.
The result [\ref{CL}] for the full cohomology groups does not hold
in our non-compact setting.
In the above example with $X=\CO$, $H^{01}\cpct(X)\ne0$
although the only fixed point $0$ is isolated.}
\end{rmk}

We now consider geometric quantization on a \ka manifold $X$ with
a holomorphic $\CM$-action satisfying Assumption~\ref{PWA}.
Recall that a pre-quantum line bundle $L$ on $(X,\om)$ is a holomorphic
line bundle whose curvature is $\frac{\om}{\ii}$.
Suppose such an $L$ exists and the $\CM$-action lifts to a holomorphic 
action on $L$.

\begin{defn}\label{QTM}{\em
The {\em quantization} of $(X,\om)$ is the virtual vector space
	\be\label{qtm}
H(X)=\bigoplus^n_{q=0}(-1)^qH^q(X,\OO(L)).
	\ee}
\end{defn}

Applying Theorem~\ref{MAIN'}.1 to the pre-quantum line bundle, we obtain
	\be
H(X)=\bigoplus_{p,q}(-1)^{p+q}E\pq_1
	\ee
as virtual representations of $\CM$, where the spaces $E\pq_1$ are given by
(\ref{e1}).

Without loss of generality, we assume that the moment map $\mu$ 
is bounded from above.
Then the $\CM$-action is plus-meromorphic.
Suppose $0$ is a regular value of $\mu$.
For simplicity, we assume that the $S^1$-action on $\mu^{-1}(0)$ is free.
Then the symplectic quotient $X_0=\mu^{-1}(0)/S^1=X^{\rm s}/\CM$ is a 
smooth \ka manifold.
We construct the symplectic cuts $(X_\pm,\om_\pm)$ as the symplectic quotients
of the $S^1$-action on $X\times\CO$, where the weights on $\CO$ are $\pm1$,
respectively [\ref{Le}].
The two cuts are \ka manifolds with holomorphic $\CM$-actions.
$X_+$ is compact and $X_-$ satisfies Assumption~\ref{PWA}.
The sets of connected components of $X^\cm_\pm$ are
$F_\pm=\zero\cup\set{\al\in F}{\mu(X^C_\al)\in\RE^\pm}$, respectively,
and $X^\cm_{\pm,0}\cong X_0$, $X^\cm_{\pm,\al}\cong X^\cm_\al$
as complex manifolds [\ref{WZ}, Lemma~4.6], which we now identify.
Let $N_0\to X_0$ be the holomorphic line bundle associate to the
circle bundle $\mu^{-1}(0)\to X_0$.
Then $\CM$ acts on the fibers of $N_0$ with weight $1$.
The holomorphic normal bundles of $X_0$ in $X_\pm$ are isomorphic to 
$N_0^{\pm1}$, respectively.
Since the action of $\CM$ lifts to $L$, the pre-quantum line bundles
$L_0\to X_0$ and $L_\pm\to X_\pm$ exist.
We have the isomorphisms $L_\pm|_{X_0}\cong L_0$ and
$L_\pm|_{X_\pm-X_0}\cong L|_{\mu^{-1}(\re^\pm)}$
(see for example [\ref{WZ}, Lemma~4.9]).

\begin{prop}\label{GLU-RED}
Under the above assumptions, we have\\
1. a gluing formula under symplectic cutting
	\be\label{glu}
\ch H(X)=\ch H(X_+)+\ch H(X_-)-\dim_\co H(X_0);
	\ee
2. that quantization commutes with reduction, i.e.,
	\be\label{red}
\dim_\co H(X)^\cm=\dim_\co H(X_0).
	\ee
\end{prop}

\proof{For $\al\in F$, let
	\be
I^\pm_\al=(-1)^{\nu^\pm_\al}\int_{X^\cm_\al}{\mathrm ch}_\cm
\left(\frac{L|_{X^{\smco^\times}_\al}\otimes\det N^\pm_\al}
{\det(1-(N^\mp_\al)^*)\otimes\det(1-N^\pm_\al)}\right){\mathrm td}(X^\cm_\al).
	\ee
Then by (\ref{ind'}), we obtain
	\bea
\ch H(X_+)\eq\sum_{\al\in F_+-\zero}I^-_\al-\int_{X_0}{\mathrm ch}_\cm\left(
\frac{L_0\otimes N^{-1}_0}{1-N_0^{-1}}\right){\mathrm td}(X_0)	\label{x+-}\\
	  \eq\sum_{\al\in F_+-\zero}I^+_\al+\int_{X_0}{\mathrm ch}_\cm
\left(\frac{L_0}{1-N_0}\right){\mathrm td}(X_0)			\label{x++}
	\eea
and
	\be\label{x--}
\ch H(X_-)=\sum_{\al\in F_--\zero}I^-_\al-\int_{X_0}{\mathrm ch}_\cm
\left(\frac{L_0\otimes N_0}{1-N_0}\right){\mathrm td}(X_0).
	\ee
1. From (\ref{x+-}) and (\ref{x--}), we get
	\be
\ch H(X_+)+\ch H(X_-)=\sum_{\al\in F}I^-_\al+
\int_{X_0}{\mathrm ch}(L_0){\mathrm td}(X_0)=\ch H(X)+\dim_\co H(X_0).
	\ee
2. From (\ref{x++}) and (\ref{x--}), we get
	\be
\dim_\co H(X_+)^T=\dim_\co H(X_-)^T=
\int_{X_0}{\mathrm ch}(L_0){\mathrm td}(X_0)=\dim_\co H(X_0).
	\ee
The result follows.}

\begin{rmk}{\em
1. We can define $H\cpct(X)=\bigoplus^n_{q=0}(-1)^qH^q\cpct(X,\OO(L))$
as the counterpart of (\ref{qtm}) with compact support.
Using (\ref{x++}) and
	\be
\ch H\cpct(X_-)=\sum_{\al\in F_--\zero}I^+_\al+\int_{X_0}{\mathrm ch}_\cm
\left(\frac{L_0}{1-N_0^*}\right){ \mathrm td}(X_0).
	\ee
we can show a similar gluing formula
	\be
\ch H\cpct(X)=\ch H(X_+)+\ch H\cpct(X_-)-\dim_\co H(X_0).
	\ee
However $\dim_\co H\cpct(X)^\cm\ne\dim_\co H(X_0)$ in general.
For example, take $X=\CO$ and choose the moment map $\mu(z)=-1-\hf|z|^2$.
Then $\dim_\co H\cpct(X)^\cm=1$ but $X_0=\emptyset$.\\
2. When $(X,\om)$ is symplectic, the individual cohomology groups in 
(\ref{qtm}) do not make sense, but $H(X)$ can be defined as the index of
a spin$^\co$-Dirac operator.
In [\ref{DGMW}], (\ref{glu}) and (\ref{red}) were proved for compact 
symplectic manifolds.
(\ref{red}) is the $S^1$-case of a conjecture by
Guillemin and Sternberg [\ref{GS}]; the cases with higher rank torus and
non-Abelian group actions were proved by Meinrenken [\ref{Me1}, \ref{Me2}],
Jeffrey and Kirwan [\ref{JK}], Vergne [\ref{V}] and others under various
generalities using localization techniques,
and by Tian and Zhang [\ref{TZ}] using an analytic approach.\\
3. Proposition~\ref{GLU-RED} shows that the results of [\ref{DGMW}] holds
for non-compact \ka manifolds under Assumption~\ref{PWA}.
For non-compact symplectic manifolds satisfying Assumption~\ref{PWA},
the validity of (\ref{glu}) and (\ref{red}) remains open.
Also it would be interesting to investigate, analytically or otherwise,
whether the Tian-Zhang inequalities [\ref{TZ}]
	\be
\sum_{k=0}^{n-1}t^k\dim_\co H^k(X_0,\OO(L_0))=
\sum_{k=0}^nt^k\dim_\co H^k(X,\OO(L))^\cm+(1+t)Q_0(t)
        \ee
for some $Q_0(t)\ge0$ hold when $X$ is non-compact and K\"ahler.
The conjecture in [\ref{WZ}, Remark 4.11] can also be posed in this
non-compact setting.}
\end{rmk}

\begin{rmk}{\em
In ordinary Morse theory, the underlying real manifold is (the bosonic 
part of) the configuration space of a supersymmetric system [\ref{W1}].
In holomorphic Morse theory, the complex manifold $X$, if it is K\"ahler,
can be interpreted as the phase space of a bosonic system;
this interpretation is adopted in Definition~\ref{QTM}.
The spectral sequence in Theorem~\ref{MAIN'}.1 or Corollary~\ref{DISCR'}.1
that converges to the quantum Hilbert space (\ref{qtm}) is a finite
dimensional model of the BRST approach in conformal field theory [\ref{Fe}].
In [\ref{FF}, \ref{BMP}], the case of flag manifolds (see subsection~4.2)
was considered.
Here we show that the analogy works for any quantizable \ka manifold with 
a Hamiltonian $S^1$-action satisfying Assumption~\ref{PWA}.
It would be interesting to extend the present work to infinite dimensional
settings.}
\end{rmk}

\vskip4ex

\noindent {\bf Acknowledgement.}
Part of this work was completed at School of Mathematics, Institute for 
Advanced Study, Princeton.
The author would like to thank M.\ Bo\v zi\v cevi\'c, M.\ Eastwood and
J.\ McCarthy for helpful discussions.

\bigskip

        \newcommand{\athr}[2]{{#1}.\ {#2}}
        \newcommand{\au}[2]{\athr{{#1}}{{#2}},}
        \newcommand{\an}[2]{\athr{{#1}}{{#2}} and}
        \newcommand{\jr}[6]{{#1}, {\it {#2}} {#3}\ ({#4}) {#5}-{#6}}
        \newcommand{\pr}[3]{{#1}, {#2} ({#3})}
        \newcommand{\bk}[4]{{\it {#1}}, {#2}, ({#3}, {#4})}
        \newcommand{\cf}[8]{{\it {#1}}, {#2}, {#5},
                 {#6}, ({#7}, {#8}), pp.\ {#3}-{#4}}
        \vspace{5ex}
        \begin{flushleft}
{\bf References}
        \end{flushleft}
{\small
        \begin{enumerate}
        
	\item\label{Ak}
	\au{E}{Akyildiz}
	\jr{Bruhat decomposition via $G_m$-action}
	{Bull.\ l'Acad.\ Polonaise Sci.}{28}{1980}{541}{547}

        \item\label{AB}
        \an{M.\ F}{Atiyah} \au{R}{Bott}
        \jr{A Lefschetz fixed point formula for elliptic complexes, Part I}
        {Ann.\ Math.}{86}{1967}{374}{407};
        \jr{Part II}{Ann.\ Math.}{87}{1968}{451}{491}

	\item\label{AS}
	\an{M.\ F}{Atiyah} \au{I.\ M}{Singer}
	\jr{The index of elliptic operators.\ III}
	{Ann.\ Math.}{87}{1968}{546}{604}

	\item\label{Au}
	\au{M}{Audin}
	\bk{The topology of torus actions on symplectic manifolds,
	{\rm Prog.\ in Math., 93}}
	{Birkh\"auser Verlag}{Basel, Boston, Berlin}{1991}

	\item\label{BS}
	\an{C}{B\u anic\u a} \au{O}{St\u an\u a\c sil\u a}
	\bk{Algebraic methods in the global theory of complex spaces}
	{Editura Academiei and John Wiley \& Sons}
	{Bucharest and London, New York, Sydney}{1976}

	\item\label{BGG}
	\au{I.\ N}{Bernstein} \an{I.\ M}{Gelfand} \au{S.\ I}{Gelfand}
	\cf{Differential operators on the base affine space and a study
	of $\sg$-modules}
	{Lie groups and their representations 
	(Bolyai J\'anos Math.\ Soc., 1971)}
	{21}{64}{ed.\ I.\ M.\ Gelfand}{Akad\'emiai Kiad\'o and Adam Hilger}
	{Budapest and London}{1975}

	\item\label{BB1}
	\au{A}{Bia\l ynicki-Birula}
	\jr{Some theorems on actions of algebraic groups}
	{Ann.\ Math.}{98}{1973}{480}{497}

	\item\label{BB2}
	\au{A}{Bia\l ynicki-Birula}
	\jr{Some properties of the decompositions of algebraic varieties 
	determined by actions of a torus}
	{Bull.\ l'Acad.\ Polonaise Sci.}{24}{1976}{667}{674}

	\item\label{B1}
	\au{R}{Bott}
	\jr{Homogeneous vector bundles}
	{Ann.\ Math.}{66}{1957}{203}{248}

	\item\label{B2}
	\au{R}{Bott}
	\jr{Morse theory indomitable}
	{I.H.E.S.\ Publ.\ Math.}{68}{1988}{99}{114}

	\item\label{BMP}
	\au{P}{Bouwknegt} \an{J}{McCarthy} \au{K}{Pilch}
	\jr{Free field realizations of WZNW models, the BRST complex and
	its quantum structure}{Phys.\ Lett.\ B}{234}{1990}{297}{303};
	\jr{Quantum group structure in the Fock space resolution of 
	$\widehat{sl}(n)$ representations}
	{Commun.\ Math.\ Phys.}{131}{1990}{125}{155};
	\cf{Free field approach to $2$-dimensional conformal field theories}
	{Common trends in mathematics and quantum field theories, Yukawa 
	International Seminar, Progr.\ Theoret.\ Phys.\ Suppl.\ No.\ 102 
	(1990)}{67}{135}{eds.\ T.\ Eguchi, T.\ Inami and T.\ Miwa} 
	{Progress of Theoretical Physics}{Kyoto}{1991}

	\item\label{Boz}
	\au{M}{Bo\v zi\v cevi\'c}
	\jr{A geometric construction of a resolution of the fundamental series}
	{Duke Math.\ J.}{60}{1990}{643}{669}

	\item\label{Br}
	\au{M}{Brion}
	\jr{Une extension du th\'eor\`eme de Borel-Weil}
	{Math.\ Ann.}{286}{1990}{655}{660}

	\item\label{B}
	\au{J.\ L}{Brylinski}
	\cf{Differential operators on the flag varieties}
	{Tableau de Young et fonctions de Schur en alg\'ebre et g\'eom\'etrie
	(Toru\'n, 1980)}{43}{60}{Ast\'erisque, v.\ 87-88}
	{Soc.\ Math.\ de France}{Paris}{1981}	

	\item\label{BK}
	\an{J.\ L}{Brylinski} \au{M}{Kashiwara}
	\jr{Kazhdan-Lusztig conjecture and holonomic systems}
	{Invent.\ Math.}{64}{1981}{387}{410}

	\item\label{CL}
	\an{J.\ B}{Carrell} \au{D}{Lieberman}
        \jr{Holomorphic vector fields and Kaehler manifolds}
        {Invent. Math.}{23}{1973}{303}{309}

	\item\label{CS1}
	\an{J.\ B}{Carrell} \au{A.\ J}{Sommese}
	\jr{{\sf C}$^*$-actions}
	{Math.\ Scad.}{43}{1978}{49}{59};
	Correction to ``\,{\sf C}$^*$-actions'', {\em idbd.} 53 (1983) 32

	\item\label{CS2}
	\an{J.\ B}{Carrell} \au{A.\ J}{Sommese}
	\jr{Some topological aspects of ${\bf C}^*$ actions on compact Kaehler
	manifolds}{Comment.\ Math.\ Helvetica}{54}{1979}{567}{582}
	
	\item\label{CS3}
	\an{J.\ B}{Carrell} \au{A.\ J}{Sommese}
	\jr{Filtration of meromorphic {\sf C}$^*$ actions on complex manifolds}
	{Math.\ Scad.}{53}{1983}{25}{31}

	\item\label{C}
	\au{C}{Chevalley}
	\cf{Sur les d\'ecomposition cellulaires des espaces $G/B$ {\em (1958)}}
	{with a foreword by A.\ Borel, 
	Algebraic groups and their generalizations: classical methods,
	Proc.\ Symp.\ Pure Math., 56 (1994), Part 1}{1}{23}
	{eds.\ W.\ J.\ Haboush and B.\ J.\ Parshall}
	{Amer.\ Math.\ Soc.}{Providence, RI}{1994}

        \item\label{DGMW}
        \au{H}{Duistermaat} \au{V}{Guillemin} \an{E}{Meinrenken} \au{S}{Wu}
        \jr{Symplectic reduction and Riemann-Roch for circle actions}
        {Math.\ Res.\ Lett.}{2}{1995}{259}{266}

	\item\label{EL}
	\an{M}{Eastwood} \au{C}{LeBrun}
	\jr{Thickening and supersymmetric extensions of complex manifolds}
	{Amer.\ J.\ Math.}{108}{1986}{1177}{1192};
	\jr{Fattening complex manifolds: curvature and Kodaira-Spencer maps}
	{J.\ Geom.\ Phys.}{8}{1992}{123}{146}
	
	\item\label{FF}
	\an{B.\ L}{Feigin} \au{E.\ V}{Frenkel}
	\jr{Affine Kac-Moody algebras and semi-infinite flag manifolds}
	{Commun.\ Math.\ Phys.}{128}{1990}{161}{189}

	\item\label{Fe}
	\au{G}{Felder}
	\jr{BRST approach to minimal models}
	{Nucl.\ Phys.\ B}{317}{1989}{215}{236};
	erratum, {\it Nucl.\ Phys.\ B} 324 (1989) 548

	\item\label{F}
	\au{T}{Frankel}
	\jr{Fixed points and torsion on \ka manifolds}
	{Ann.\ Math.}{70}{1959}{1}{8}

	\item\label{Fu}
	\au{A}{Fujiki}
	\jr{Fixed points of the actions on compact \ka manifolds}
	{Publ.\ R.I.M.S., Kyoto Univ.}{15}{1979}{797}{826}

	\item\label{G}
	\au{R}{Godement}
	\bk{Topologie alg\'ebrique et th\'eorie des faisceaux}
	{Hermann}{Paris}{1958}

	\item\label{Gr}
	\au{Ph.\ A}{Griffiths}
	\jr{The extension problem in complex analysis II: embeddings with
	positive normal bundle}
	{Amer.\ J.\ Math.}{88}{1966}{366}{446}

	\item\label{GS}
	\an{V}{Guillemin} \au{S}{Sternberg}
	\jr{Geometric quantization and multiplicities of group representations}
	{Invent. Math.}{67}{1982}{515}{538}

	\item\label{H}
	\au{R}{Hartshorne}
	\bk{Residues and duality, {\em Lecture Notes in Math., No.\ 20}}
	{Springer-Verlag}{Berlin}{1966}

        \item\label{JK}
        \an{L.\ C}{Jeffrey} \au{F.\ C}{Kirwan}
        \jr{Localization and the quantization conjecture}
        {Topology}{36}{1997}{647}{693}

	\item\label{J}
	\au{J}{Jurkiewicz}
	\jr{An example of algebraic torus action which determines the
	nonfiltrable decomposition}
	{Bull.\ l'Acad.\ Polonaise Sci.}{25}{1977}{1089}{1092}

	\item\label{K}
	\au{G}{Kempf}
	\jr{The Grothendieck-Cousin complex of an induced representation}
	{Adv.\ in Math.}{29}{1978}{310}{396}

	\item\label{Kon}
	\au{J}{Konarski}
	\jr{Decompositions of normal algebraic varieties determined by an
	action of a one-dimensional torus}
	{Bull.\ l'Acad.\ Polonaise Sci.}{26}{1978}{295}{300}

	\item\label{Kor}
	\au{M}{Koras}
	\jr{On actions of an analytic torus}
	{Bull.\ l'Acad.\ Polonaise Sci.}{27}{1979}{21}{26}

	\item\label{La}
	\au{M}{Landucci}
	\jr{Solutions with ``precise'' compact support of the $\pdb$-problem 
	in strictly pseudoconvex domains and some consequences}
	{Atti Accad.\ Naz.\ Lincei Rend.\ Cl.\ Sci.\ Fis.\ Mat.\ Natur.\ (8)}
	{67}{1979}{81}{86};
	\jr{Solutions with precise compact support of $\pdb u=f$}
	{Bull.\ Sci.\ Math.\ (2)}{104}{1980}{273}{299}

	\item\label{Lep}
	\au{J}{Lepowsky}
	\jr{A generalization of the Bernstein-Gelfand-Gelfand resolution}
	{J.\ Algebra}{49}{1977}{496}{511}

	\item\label{Le}
	\au{E}{Lerman}
	\jr{Symplectic cuts}{Math.\ Res.\ Lett.}{2}{1995}{247}{258}

	\item\label{MW}
	\an{V}{Mathai} \au{S}{Wu}
	\jr{Equivariant holomorphic Morse inequalities I: a heat kernel proof}
	{J.\ Diff.\ Geom.}{46}{1997}{78}{98}

        \item\label{Me1}
        \au{E}{Meinrenken}
        \jr{On Riemann-Roch formulas for multiplicities}
        {J.\ Amer.\ Math.\ Soc.}{9}{1996}{373}{389}

        \item\label{Me2}
        \au{E}{Meinrenken}
        \jr{Symplectic surgery and the spin$^c$-Dirac operator}
        {Adv.\ Math.}{134}{1998}{240}{277}

	\item\label{MR}
	\an{M}{Murray} \au{J}{Rice}
	\jr{A geometric realisation of the Lepowsky Bernstein Gelfand Gelfand
	resolution}
	{Proc.\ Amer.\ Math.\ Soc.}{114}{1992}{553}{559}

	\item\label{P}
	\au{E}{Prato} 
	\jr{Convexity properties of the moment map for certain non-compact
	manifolds}{Comm.\ Anal.\ Geom.}{2}{1994}{267}{278}

	\item\label{PW}
	\an{E}{Prato} \au{S}{Wu}
	\jr{Duistermaat-Heckman measures in a non-compact setting}
	{Comp.\ Math.}{94}{1994}{113}{128}
	
	\item\label{Sm}
	\au{S}{Smale}
	\jr{Differential dynamical systems}
	{Bull.\ Am.\ Math.\ Soc.}{73}{1967}{747}{817}

	\item\label{So}
	\au{A.\ J}{Sommese}
	\cf{Some examples of $\,C^*\!$ actions}
	{Group actions and vector fields (Proc.\ of Polish-North American 
	Seminar, Univ.\ of British Columbia, 1981), 
	Lecture Notes in Math., 956}{118}{124}{ed.\ J.\ B.\ Carrell}
	{Springer-Verlag}{Berlin, Heidelberg, New York}{1982}

        \item\label{TZ}
        \an{Y}{Tian} \au{W}{Zhang}
        \jr{Symplectic reduction and quantization}
        {C. R. Acad. Sci. Paris, S\'erie I}{324}{1997}{433}{438};
        \jr{An analytic proof of the geometric quantization conjecture of
        Guillemin-Sternberg}{Invent.\ Math.}{132}{1998}{229}{259}

        \item\label{V}
        \au{M}{Vergne}
        \jr{Multiplicities formula for geometric quantization. I, II}
        {Duke Math. J.}{82}{1996}{143}{179}, 181-194

	\item\label{Wa}
	\au{G}{Warner}
	\bk{Harmonic analysis on semi-simple Lie groups I,
	{\rm Grund.\ Math.\ Wiss.\ 188}}
	{Springer-Verlag}{Berlin, Heidelberg, New York}{1992}

	\item\label{W1}
	\au{E}{Witten} 
	\jr{Supersymmetry and Morse theory}
	{J.\ Diff.\ Geom.}{17}{1982}{661}{692}
	
	\item\label{W2} 
	\au{E}{Witten} 
	\cf{Holomorphic Morse inequalities}
	{Algebraic and differential topology, Teubner-Texte Math., 70}
	{318}{333}{ed.\ G.\ Rassias}{Teubner}{Leipzig}{1984}

	\item\label{Wu}
	\au{S}{Wu}
	\pr{Equivariant holomorphic Morse inequalities II: torus and
	non-Abelian group actions}
	{MSRI preprint No.~1996-013, {\tt dg-ga/9602008}}{1996}

	\item\label{WZ}
	\an{S}{Wu} \au{W}{Zhang}
	\jr{Equivariant holomorphic Morse inequalities III: non-isolated 
	fixed points}{Geom.\ Funct.\ Anal.}{8}{1998}{149}{178}

	\item\label{Z}
	\au{G.\ J}{Zuckerman}
	\cf{Geometric methods in representation theory}
	{Representation theory of reductive groups, Prog.\ in Math., 40}
	{283}{290}{ed.\ P.\ C.\ Trombi}{Birkh\"auser}{Boston}{1983}
	
        \end{enumerate}}
	\end{document}